\documentclass{article}

\usepackage{arxiv}

\usepackage{amsmath}
\usepackage{graphicx}
\usepackage{epstopdf}
\usepackage[utf8]{inputenc} % allow utf-8 input
\usepackage[T1]{fontenc}    % use 8-bit T1 fonts
\usepackage{hyperref}       % hyperlinks
\usepackage{url}            % simple URL typesetting
\usepackage{booktabs}       % professional-quality tables
\usepackage{amsfonts}       % blackboard math symbols
\usepackage{nicefrac}       % compact symbols for 1/2, etc.
\usepackage{microtype}      % microtypography
\usepackage{lipsum}
\usepackage{algorithm}
\usepackage{algorithmic}
\usepackage{booktabs}
\usepackage{lscape}
\usepackage{caption}
\usepackage{subcaption}
\usepackage{multirow}
\newtheorem{thm}{Theorem}[section]
\newtheorem{lem}[thm]{Lemma}

\newtheorem{defi}{Definition}[section]

\newtheorem{claim}{Claim}[section]

\newcommand{\goto}{\ensuremath{\rightarrow}}

\newenvironment{Assumptions}
{
\setcounter{enumi}{0}

\begin{enumerate}}
{\end{enumerate} }

\numberwithin{equation}{section} \allowdisplaybreaks

\title{A Gray Level Indicator-Based Regularized Telegraph Diffusion Equation Applied to Image Despeckling}

\author{
 Sudeb Majee \\
  School of Basic Sciences\\ 
  Indian Institute of Technology Mandi\\
  PIN 175005, INDIA\\
  \texttt{sudebmajee@gmail.com} \\
    \And
  Rajendra K. Ray \\
  School of Basic Sciences\\ 
  Indian Institute of Technology Mandi\\
  PIN 175005, INDIA\\
  \texttt{rajendra@iitmandi.ac.in} \\
    \And
   Ananta K. Majee \\
   Department of Mathematics\\
   Indian Institute of Technology Delhi\\
   PIN 110016, INDIA \\
  \texttt{majee@maths.iitd.ac.in} \\ 
}

\begin{document}
\maketitle

\begin{abstract}
In this work, a gray level indicator based non-linear telegraph diffusion model is presented for multiplicative noise removal problem. Most of the researchers focus only on diffusion equation-based model for multiplicative noise removal problem. The suggested model uses the benefit of the combined effect of diffusion equation as well as the wave equation. Wave nature of the model preserves the high oscillatory and texture pattern in an image. In this model, the diffusion coefficient depends not only on the image gradient but also on the gray level of the image, which controls the diffusion process better than only gradient-based diffusion models.
Moreover, we prove the well-posedness of the present model using Schauder fixed point theorem. Furthermore, we show the superiority of the proposed model over a recently developed method on a set of gray level test images which are corrupted by speckle noise.
\end{abstract}

\section{Introduction}
\label{intro}
In the real scenario, images are often corrupted by different types of noises, e.g., additive, multiplicative, or mixed nature. Hence the noise removal process is a very initial stage for high-level image
analysis. In this work, we focus our interest only on multiplicative speckle noise removal process. Purity of the edge/texture information in the synthetic aperture radar(SAR) images, ultrasound images, and
laser images are usually diminished by speckle noise \cite{burckhardt1978speckle,loizou2005comparative,prager2001speckle}. Due to the contamination by speckle noise, it is challenging to distinguish
the hidden details in the images. Therefore the development of an advance speckle noise removal algorithm is always an essential aspect for the image processing society. A Mathematical representation
is still required to develop an efficient noise removal algorithm so that we can express each pixel of an image as a function of the speckle noise. The popularly used model for the noise image can be express
as a product of the original signal and the speckle noise \cite{dutt1995statistical}
\begin{equation*}
J=I\eta,
\end{equation*}
where $J$ indicates the noisy image, $I$ is the noise-free image, and $\eta$  signifies the speckle-noise process.

In general, the probability density function of the multiplicative speckle noise process $\eta$  follows the Gamma Law as,
\begin{equation*}
g(\eta) =
        \begin{cases}
        \frac{L^L}{\Gamma \left( L \right)}\eta^{L-1}\text{exp}\left(-L\eta\right), \hspace{0.5cm} \text{for}\hspace{0.1cm} \eta >0, \\
        \hspace{1.5cm}  0, \hspace{1.9cm} \text{for}\hspace{0.1cm} \eta = 0,
        \end{cases}
\end{equation*}
where $L \in {\rm I\!N} $ signifies the number of looks which correspond to the noise level in the corrupted images \cite{argenti2013tutorial,hao2015variational,liu2016modified} and $\Gamma \left( \cdot \right)$ denotes the Gamma function.

A large number of study reports the fundamentals and the statistical attributes of the speckle noise \cite{achim2001novel,argenti2013tutorial,frost1982model,kuan1985adaptive,lee1980digital,prager2001speckle}.
Futuristic despeckling approaches include Bayesian methods in spatial domain \cite{frost1982model,kuan1985adaptive,lee1980digital},
Bayesian methods in transformed domain \cite{aiazzi1998multiresolution,hao1999novel,meer1994multiresolution}, order-statistics and
morphological filters \cite{alparone1998decimated,alparone1995two,crimmins1985geometric,prager2001speckle},
simulated annealing despeckling \cite{white1994simulated}, nonlocal
filtering \cite{buades2005non,coupe2008bayesian,deledalle2009iterative,teuber2012new},
wavelet-based approaches \cite{achim2001novel,sudha2009speckle}, 
nonlinear diffusion in Laplacian pyramid domain \cite{zhang2007nonlinear},
anisotropic diffusion based methods \cite{jain2019non,jain2018nonlinear,jin2000adaptive,shan2019multiplicative,yu2002speckle,zhou2015doubly,zhou2018nonlinear},
and variational methods \cite{  aubert2008variational,dong2013convex,huang2010multiplicative,jidesh2013complex, jin2010analysis, jin2011variational, liu2013nondivergence,rudin2003multiplicative,shi2008nonlinear}.

From the initiation of PM model \cite{perona1990scale}, the partial differential equations(PDEs) are extensively used to develop noise removal algorithms, among different types of PDE based models, the total variational (TV)  based algorithms are achieved remarkable results. First variational based strategy to deal with multiplicative noise is proposed by Rudin et al. \cite{rudin2003multiplicative}, with the principles,  
\begin{equation*}
\int_{\Omega} \frac{J}{I}dx=1, \quad {\rm and}\quad  \int_{\Omega} {\left( \frac{J}{I}-1\right) }^2dx=\sigma^2\,,
\end{equation*}
where $\sigma^2$ represents the variance of the noise $\eta$.
Due to the non-convexity of their proposed energy function, the model may not give a globally unique solution. To succeed over this shortcoming, several authors suggested various convex functional with different data fidelity terms \cite{aubert2008variational,huang2010multiplicative, jin2011variational, liu2013nondivergence}.
Recently, Dong et al. \cite{dong2013convex} suggest a convex total variation model for multiplicative speckle-noise reduction with the following form:
\begin{equation}\label{eq:Dong_energy}
I=\min_{I\in \text{BV} \left(\Omega\right)}\left\lbrace \int_{\Omega} \alpha(x)|\nabla I|dx+\lambda \int_{\Omega}\left( I+J \log\frac{1}{I}\right) dx\right\rbrace. \nonumber
\end{equation} 
They choose the gray level indicator function $\alpha$, as
\begin{equation}\label{eq:gray_indicator}
\left( 1-\frac{1}{1+k|G_\xi \ast J|^2}\right) \frac{1+kM^2}{kM^2},\,\,\,\, \text{or}\,\,\,\,\, \dfrac{G_{\xi}\ast J}{M}, \nonumber
\end{equation}
with $M=\underset{x \in \Omega}{\text{sup}}(G_\xi \ast J)(x)$, where $\xi>0$, $k>0$, ``$\ast$" is the convolution operator, $G_\xi$ is the two dimensional Gaussian kernel and $\lambda$ is a given parameter, see \cite{dong2013convex}. Later, based on a gray level indicator function, Zhou et al.  proposed a diffusion model~(DDD model)\cite{zhou2015doubly} for multiplicative noise removal problem. Their model takes the form:
\begin{align*}
&I_t  = \text{div}(g(I,|\nabla I|)\nabla I),  \hspace{0.3cm} \text{in} \hspace{0.2cm} \Omega_T:= \Omega \times (0,T),  \\
&\partial_n I=0,              \hspace{2.6cm} \text{in} \hspace{0.2cm} \partial \Omega_T:= \partial \Omega \times (0,T),\\
&I(x,0)=I_0(x),	\hspace{1.5cm} \text{in} \hspace{0.2cm}\Omega,
\end{align*}
where $\Omega$ is the domain of original image $I$ and the observed noise image $I_0$, $\text{div}$ and $\nabla$ represents the divergence and gradient operator respectively.
They choose the diffusion coefficient as
\begin{align*}
g\left(I,\vert \nabla I \vert \right)=\dfrac{2\vert I \vert^\nu}{M^\nu+\vert I \vert^\nu}\cdot \dfrac{1}{\left(1+ |\nabla I|^2\right)^{(1-\beta)/2} }, 
\end{align*}
where $\nu>0,$ $0<\beta<1,$ and $M=\underset{x \in \Omega}{\text{sup}} I$. In this case, the gray level indicator and edge detector function are $a(I):=\dfrac{2\vert I \vert^\nu}{M^\nu+\vert I \vert^\nu}$  and  $b(I):= \dfrac{1}{\left(1+ |\nabla I|^2\right)^{(1-\beta)/2} }$ respectively.  
However because of the degeneracy of the edge detector function, i.e., $b(|\nabla I |) \rightarrow 0$ as $|\nabla I | \rightarrow \infty$, it is challenging to establish the well-posedness of their model.
Recently Shan et al. \cite{shan2019multiplicative} proposed a regularized version of the above-discussed model \cite{zhou2015doubly}.
In \cite{shan2019multiplicative}, the model takes of the form
\begin{align*}
&I_t  = \text{div}(g(I_\xi,|\nabla I_\xi|)\nabla I),  \hspace{0.2cm} \text{in} \hspace{0.2cm} \Omega_T,  \\
&\partial_n I=0,              \hspace{2.6cm} \text{in} \hspace{0.2cm} \partial \Omega_T,\\
&I(x,0)=I_0(x),	\hspace{1.5cm} \text{in} \hspace{0.2cm}\Omega\,.
\end{align*}
They choose the diffusion coefficient as
\begin{align*}
g\left(I,\vert \nabla I \vert \right)= \left( \dfrac{I_\xi}{M_\xi^I} \right)^\nu \cdot \dfrac{1}{1+ |\nabla I_\xi|^\beta }, 
\end{align*}
where $I_\xi=G_\xi\ast I$, $M_\xi^I= \underset{x\in \Omega}{\text{max}}\vert I_\xi(x,t) \vert$, $I_0$ is  the initial image, and $\nu, \beta$ and $\xi$ are positive constants. Due to the introduction of the Gaussian kernel in the diffusion coefficient, which avoids the degeneracy in the model, the authors able to study the wellposedness of the underlying problem.

To the best of our knowledge, most of the researcher concentrated their interest only on parabolic PDE based models, which are developed either from the variational bases approach or diffusion based approach, for the speckle-noise removal process. The hyperbolic PDEs could upgrade the quality of the detected edges and improve the image better than parabolic PDEs \cite{averbuch2006edge}. In the existing literature, the first
hyperbolic model for image denoising is telegraph-diffusion model \cite{ratner2007image}, where the image was viewed as an elastic sheet placed in a damping environment, which interpolates between the diffusion
equation and the wave equation. The telegraph-diffusion model takes the form,
\begin{align*}
&I_{tt}+\gamma I_t  =\text{div}(g(|\nabla I|)\nabla I),  \hspace{0.6cm} \text{in} \hspace{0.2cm} \Omega_T,  \\
&\partial_n I=0,      \hspace{3.5cm} \text{in}
\hspace{0.2cm} \partial \Omega_T,\\
&I(x,0)=I_0(x), \hspace{0.2cm} I_t(x,0)=0,	\hspace{0.3cm} \text{in} \hspace{0.2cm}\Omega,
\end{align*}
where $g(|\nabla I|)=1/(1+({|\nabla I|^2}/{k^2}))$ is an edge-controlled diffusion function which preserves the important features and smoothen the unwanted signals, and $\gamma$ is the damping parameter.
It is quite interesting to note that for a very higher value of $g$ and $\gamma$, this telegraph-diffusion equation (TDE model) converges
to the original PM model \cite{perona1990scale} in a long time scenario. Although the TDE model performs better, it is challenging to confirm the well-posedness of their model. To overcome the ill-posedness issue in the TDE model \cite{ratner2007image}, Cao et al. suggest a regularized TDE model \cite{cao2010class}. They replace the gradient $|\nabla I|$  by $|\nabla G_\xi \ast I| $ in the edge-controlled function $g$ in the TDE model \cite{ratner2007image} and establish the well-posedness of their proposed model. 
Even though the TDE model can effectively preserve the sharp edges but failed to produce satisfactory smoothing in the presence of a large level of noise.
To overcome this issue, several non-linear telegraph diffusion-based method have been proposed \cite{cao2010class,jain2016edge,sun2016class,yang2014kernel,zhang2015spatial}. However, in spite of their impressive applications in the additive noise removal process, hyperbolic PDE based approaches have not successfully used for speckle noise removal process.

Recently Sudeb et al. suggest a fuzzy edge detector based telegraph total variation model \cite{fuzzy2019ttvmodel} for the speckle noise removal problem. To the best of our knowledge, this is the first hyperbolic PDE based model in the existing literature applied to speckle noise removal process. The model \cite{fuzzy2019ttvmodel} takes the form,
\begin{align*}
&I_{tt}+\gamma I_t = \text{div}\left(\theta(I)\frac{\nabla I}{|\nabla I|}\right)-\lambda \left( 1-\frac{I_0}{I}\right),
\hspace{0.6cm} \text{in} \hspace{0.2cm} \Omega_{T},  \\
&\partial_n I =0,               \hspace{6.1cm} \text{in}
\hspace{0.2cm} \partial \Omega_{T},\\
&I(x,0)=I_0(x), \hspace{0.1cm}   I_t(x,0)=0,	\hspace{3.0cm} \text{in} \hspace{0.2cm}\Omega,
\end{align*}
where $\theta$ is the fuzzy edge detector function \cite{chaira2008new}, $\gamma$ is a positive parameter and $\lambda$ is the weight parameter.
%\begin{equation*}
%\lambda = \frac{1}{\sigma^2 |\Omega|} \int_{\Omega} \text{div} \left(\theta(I)\frac{\nabla I}{|\nabla I|}\right)\left( 1-\frac{I_0}{I}\right)dx.
%\end{equation*}

Further continuing to demonstrate the importance of hyperbolic PDE based model for image despeckling, the present work suggests a gray level indicator based telegraph diffusion model for multiplicative speckle noise removal. In this model, we choose a different diffusivity function from our previous model \cite{fuzzy2019ttvmodel}.
Also, instead of total variation framework \cite{fuzzy2019ttvmodel}, we designed the present model in an anisotropic diffusion-based fashion as discussed in \cite{zhou2015doubly}. Furthermore, we study the well-posedness of the suggested model in an appropriate function space. We opt an explicit numerical method to solve the present model. Our numerical implementation allows computing despeckled results on some standard test images. Quality of the despeckled images using the suggested model compare with the recently developed model \cite{shan2019multiplicative}.  We compare the quantitative and qualitative results at different noise levels. The experiment results  confirm that the proposed model performs better as compared to the model considered for the comparison.

The rest of the paper is organized as follows. Section \ref{sec:Proposed Model} describes the proposed telegraph diffusion method for image despeckling. In section \ref{sec:analysis}, we study the wellposedness of  weak solution of the proposed model. Section \ref{sec:numerical} describes the numerical discretization of the present model. The simulated despeckling results obtained by the proposed approach are compared with other discussed diffusion methods in Section \ref{sec:Results}. We conclude the paper in Section \ref{sec:Conclusion} with a scope on future work.

\section{Telegraph Diffusion Model for Speckle Noise Removal}
\label{sec:Proposed Model}
Inspired by the ideas of \cite{fuzzy2019ttvmodel} and \cite{zhou2015doubly} initially we developed the model
\begin{align}
&I_{tt} +\gamma I_{t}- \text{div} \left(   g\left(I,\vert \nabla I \vert \right)    \nabla I\right)=-\lambda h(I_0,I)\,, \hspace{0.2cm} \text{in}\,\,\, \Omega_T\,, \label{maina1} \\
\label{mainb1}
&\partial_n I=0\,, \hspace{5.7cm} \text{on}\,\,\,  \partial\Omega_T\,,\\
&I(x,0)=I_0(x)\,, \hspace{0.2cm} I_t(x,0)=0\,, \hspace{2.5cm} \text{in}\,\,\, \Omega\,.\label{mainc1}
\end{align}
%with the boundary and initial conditions  \eqref{tdm:bc} and \eqref{tdm:inc} respectively.
The function $g$ is defined as
\begin{align}\label{g_ddd}
g\left(I,\vert \nabla I \vert \right)=\dfrac{2\vert I \vert^\nu}{\big( M^{I}\big)^\nu+\vert I \vert^\nu}. \dfrac{1}{1+ \left(\frac{|\nabla I|}{K} \right)^2 }, 
\end{align}
where, $\nu \geq 1,$ $\gamma, K>0$ are constants, $M^{I}= \underset{x \in \Omega}{\text{max}}\vert I(x,t) \vert$, and 
$h(I_0, I)$ is the source term which comes due to the fidelity control term in the energy functional as discussed in \cite{fuzzy2019ttvmodel}. Although the presence of fidelity term in the equation keeps the restored image close to the original image, the noise may not be removed sufficiently. Therefore we would like to choose $h(I_0, I)=0.$ Also, because of the degeneracy in the diffusion coefficient \ref{g_ddd}, the suggested model \eqref{maina1}-\eqref{mainc1} may not be a well-posed problem\cite{shan2019multiplicative}. To overcome these issues, we invoke the ideas of \cite{cao2010class} and \cite{shan2019multiplicative}, and finally designe the following model in the anisotropic diffusion-based framework:

\begin{align}
&I_{tt} +\gamma I_{t}- \text{div} \left( g\left(I_\xi,\vert \nabla I_\xi \vert \right)    \nabla I\right)=0\,, \hspace{1.1cm} \text{in}\,\,\, \Omega_T\,, \label{maina} \\
\label{mainb}
&\partial_n I=0\,, \hspace{5.4cm} \text{on}\,\,\,  \partial\Omega_T\,,\\
&I(x,0)=I_0(x)\,, \hspace{0.2cm} I_t(x,0)=0\,, \hspace{2.2cm} \text{in}\,\,\, \Omega\,,\label{mainc}
\end{align}
where the diffusion function $g$ as given by
\begin{align*}
g\left(I_\xi,\vert \nabla I_\xi \vert \right)=\dfrac{ 2\vert I_\xi \vert^\nu}{\big(M^{I}_{\xi}\big)^\nu+\vert I_\xi \vert^\nu}\cdot \dfrac{1}{1+ \left(\frac{|\nabla I_{\xi}|}{K} \right)^2 }\,.
\end{align*}
In the above, $I_\xi=G_\xi\ast I$, $M^{I}_{\xi}= \underset{ x \in \Omega}{\text{max}}\vert I_\xi(x,t) \vert.$ Moreover the gray level indicator function 
\begin{align*}
b(I)=\dfrac{2\vert I_\xi \vert^\nu}{\big(M^{I}_{\xi}\big)^\nu+\vert I_\xi \vert^\nu}
\end{align*}
can be transformed into  $ b(s)=\dfrac{2s^\nu}{1+s^\nu} $, where  $s=\dfrac{|I_{\xi}|}{M_\xi^I} \in [0, 1].$

 The use of Gaussian convolution in the proposed model has a lot of advantages, not only the robustness in denoising viewpoint but also the well-posedness in the theoretical perspective.  There are two key advantages of this proposed approach:
 \begin{itemize}
 \item[i)] it provides the sharp and true edges during noise removal process than other non-telegraph based algorithms as the model \eqref{maina}-\eqref{mainc}  consists of telegraph diffusion model \cite{ratner2007image}
 \item[ii)]  it controls the diffusion process very well along with the gradient based edge detector coefficient specially for the speckle noise removal process \cite{dong2013convex} as  the gray level indicator function in the proposed model is incorporated into the telegraph diffusion framework.
\end{itemize}

\section{Wellposedness of weak solution}
\label{sec:analysis}
In this section, we prove the existence and uniqueness of weak solution of the proposed model \eqref{maina}-\eqref{mainc}. Since
the problem \eqref{maina}-\eqref{mainc} is nonlinear, we first consider the linearized problem, and then use Schauder's
fixed-point theorem  \cite{LCEvans1998} to show the existence of a weak solution. Without loss of generality, we assume $\gamma=1$ in \eqref{maina}.
%%%%%%%%%%%%%%%%%%%%%%%%%%%%%%%%%%%%%%%%%%%%%%%%%%%%%%%%%%%%%%%%%%%%%%%%%%%%%%%%%%%%%%%%%%%%%%%%%%%%%%%%%%%%%%%%%%%%%%%%%%%%%%%%%%%%%%%%%%%%%%%%%%%%%%%%%%%%%%%%%%%%%%%%%%%%%%%%%%%%%%%%%%%%%%%%%%%%%%
\subsection{Technical framework $\&$ statement of the main result}
Throughout this section, $C$ denotes a generic positive constant. For $1\le p\le \infty$, we denote by $(L^p, \|\cdot\|_{L^p})$ the standard spaces of $p$-th order integrable
functions on $\Omega$. For $r\in \mathbb{N}$,
we write $(H^r, \|\cdot\|_{H^r})$ for usual Sobolev spaces on $\Omega$, and $(H^{1})^\prime$ for the dual space of $H^1$. 
%%%%%%%%%%%%%%%%%%%%%%%%%%%%%%%%%%%%%%%%%%%%%%%%%%%%%%%%%%%%%%%%%%%%%%%%%%%%%%%%%%%%%%%%%%%%%%%%%%%%%%%%%%%%%%%%%%%%%%%%%%%%%%%%%%%%%%%
We introduce the solution space $W(0,T)$ for the
problem \eqref{maina}-\eqref{mainc}, where
\begin{align*}
 W(0,T)&=\Big\{w\in L^\infty(0,T; H^1)\,, w_t \in L^\infty(0,T; L^2); \,w_{tt} \in L^2(0,T; (H^1)') \Big\}\,.
\end{align*}
Note that the space $W(0,T)$ is a Hilbert space for the graph norm, see \cite{jllions1968}. 
\begin{defi}[Weak solution]\label{defi:weak}
A function $I$ is called a weak solution of \eqref{maina}-\eqref{mainc} if
\begin{itemize}
 \item[a)] $I \in W(0,T) $ and \eqref{mainc} holds.
 \item[b)] For all $\phi \in H^1$ and a.e $t\in (0,T)$, there hold
 \begin{align*}
\left \langle I_{tt}, \phi   \right \rangle + {\displaystyle \int_{\Omega}}
\Big(  I_t \phi +  g\left(I_\xi,\vert \nabla I_\xi \vert \right)  \nabla I\cdot \nabla \phi  \Big)\,dx =0.
  \end{align*}
\end{itemize}
\end{defi}
As we mentioned, our aim is to establish wellposedness of weak solutions of the underlying problem \eqref{maina}-\eqref{mainc}, and we will do so under the following assumption:
\begin{Assumptions}
 \item \label{A1}  The initial data $I_0$ is an $H^2$-valued function such that
 \begin{align*}
  0< \alpha:=\inf_{x\in \Omega} I_0(x)\,.
 \end{align*}
\end{Assumptions}
\begin{thm}\label{thm:existence-uniqueness}
Let the assumption \ref{A1} be true. Then the problem \eqref{maina}-\eqref{mainc} admits  a unique weak solution in the sense of Definition \ref{defi:weak}. 
\end{thm}
%%%%%%%%%%%%%%%%%%%%%%%%%%%%%%%%%%%%%%%%%%%%%%%%%%%%%%%%%%%%%%%%%%%%%%%%%%%%%%%%%%%%%%%%%%%%%%%%%%%%%%%%%%%%%%%%%%%%%%%%%%%%%%%%%%%%%%%%%%%%%%%%%%%%%%%%%%%%%%
\subsection{Linearized problem  $\&$ existence of weak solution:} For any positive constant $M_1>0$, define
\begin{align*}
 W_{M_1}= \big\{ & \bar{I}\in W(0,T):~  \|\bar{I}\|_{L^\infty(0,T;H^1)} + \|\bar{I}_t\|_{L^\infty(0,T; L^2)} \le M_1\|I_0\|_{H^1},\,\\
 & \hspace{0.5cm} 0<\alpha \le \bar{I}(x,t)~~{\rm  for~a.e.}~~(x,t)\in \Omega_T\,
  \big\}.  
\end{align*}
For any $\bar{I}\in W_{M_1}$, consider the linearized problem:
\begin{align}
I_{tt} + I_{t}- {\rm div}\big( \bar{g}(x,t) \nabla I\big)=0 \hspace{1cm} {\rm in}~~~ \Omega_T\,, \label{linmaina}
\end{align}
with the initial condition \eqref{mainc}, 
where the function $\bar{g}$ is given by
\begin{align*}
 \bar{g}(x,t) \equiv g_{\bar{I}}(x,t):= \frac{|\bar{I}_\xi|^\nu}{ \big(M_\xi^{\bar{I}}\big)^\nu + |\bar{I}_\xi|^\nu}\cdot \dfrac{1}{1+ \left(\frac{|\nabla \bar{I}_{\xi}|}{K} \right)^2 }\,.
\end{align*}
\begin{claim}\label{claim:1}
 There exist positive constants $\kappa, C >0$, depending only on $G_\xi, I_0, M_1, K, \alpha$ and $\nu$, such that
\begin{equation}\label{bound:g_w}
 \begin{aligned}
  &{\rm i)}~ 0< \kappa \le \bar{g}\le 1\,, \\
  & {\rm ii)}~ |\bar{g}_t| \le C\,.
 \end{aligned}
\end{equation}
\end{claim}

\textbf{Proof:}
%\begin{proof}
 \noindent{Proof of ${\rm i)}$:} Since $\bar{I} \in W_{M_1}$, by convolution property, we have
 \begin{align*}
  \alpha \|G_\xi\|_{L^1} \le | G_\xi \ast \alpha| \le |\bar{I}_\xi| \le M_1 C_\xi \|I_0\|_{H^1}\,; \quad 
  \big(\alpha\|G_\xi\|_{L^1}\big)^\nu \le \big(M_\xi^{\bar{I}}\big)^\nu \le \big( M_1 C_\xi \|I_0\|_{H^1}\big)^\nu\,,
 \end{align*}
and hence
\begin{align}
 \frac{ \big(\alpha \|G_\xi\|_{L^1}\big)^\nu}{2\,\big( M_1 C_\xi \|I_0\|_{H^1}\big)^\nu} \le \frac{|\bar{I}_\xi|^\nu}{ \big(M_\xi^{\bar{I}}\big)^\nu + |\bar{I}_\xi|^\nu} \le 1\,. \label{esti:1-bar-g}
\end{align}
Again by Young's convolution inequality, we observe that
\begin{align}
\dfrac{1}{1+\left( \frac{C_{\xi} M_1\left\Vert I_0 \right\Vert_{H^1}}{K} \right)^2 } \leq \frac{1}{1+ \left(\frac{|\nabla \bar{I}_{\xi}|}{K} \right)^2 }\leq 1. \label{esti:2-bar-g}
\end{align}
Now ${\rm i)}$ follows from \eqref{esti:1-bar-g}-\eqref{esti:2-bar-g} for $\kappa=\frac{ \big(\alpha \|G_\xi\|_{L^1}\big)^\nu}{2\,\big( M_1 C_\xi \|I_0\|_{H^1}\big)^\nu}
\cdot \dfrac{1}{1+\left( \frac{C_{\xi} M_1\left\Vert I_0 \right\Vert_{H^1}}{K} \right)^2 }   $. 
\vspace{.1cm}

\noindent{Proof of ${\rm ii)}:$} Observe that, since $0< \alpha \|G_\xi\|_{L^1} < M_\xi^{\bar{I}}$, we have
\begin{align*}
|\bar{g}_t| & \le C(\nu,\alpha,\xi, M_1,\Vert I_0 \Vert_{H^1}) + C(\xi, K, M_1)\Vert I_0 \Vert^2_{H^1}   \,.
\end{align*}
Thus ${\rm ii)}$ holds. This finishes the proof of claim. 
%\end{proof}
%%%%%%%%%%%%%%%%%%%%%%%%%%%%%%%%
Thanks to Claim \ref{claim:1},  one can apply classical Galerkin method \cite{LCEvans1998} to show that there exists a 
unique weak solution $ I \in W(0,T)$ of the  linearized problem \eqref{linmaina} with the initial condition \eqref{mainc}. 
\begin{lem}\label{lem:a-priori}
The unique solution $ I \in W(0,T)$ of  the  linearized problem \eqref{linmaina} with the initial condition \eqref{mainc} satisfies the following: there exists a constant $C>0$, depending only on 
$G_\xi, I_0, M_1, \nu, \alpha, K$ such that
\begin{itemize}
\item[a)] $ \|I\|_{L^\infty(0,T; H^1)} + \|I_t\|_{L^\infty(0,T; L^2)} \le C \|I_0\|_{H^1}$, \\
\item[b)] $\int_0^T \|I_{tt}\|_{(H^1)^\prime}^2\,dt \le C T \|I_0\|_{H^1}^2$. \\
\end{itemize}
\end{lem}
%=======================================================================================================================
%\begin{proof}
\textbf{Proof:}
\noindent{Proof of ${\rm a)}$:}
Note that $I_t \in L^\infty(0,T; H^1)$. Taking $\phi=I_t$ in \eqref{linmaina}, integrating by parts and using the 
inequality $\int_{\Omega} \bar{g}\nabla I \cdot \nabla I_t\, dx 
\geq  \frac{1}{2}\dfrac{d}{dt}\int_{\Omega} \bar{g}|\nabla |^2 \,dx- \frac{C}{2}\|\nabla I\|_{L^2}^2$, 
which follows from 
integration by parts formula and \eqref{bound:g_w}, and the fact 
\begin{align}
\|\nabla I\|_{L^2}^2 \le \frac{1}{\kappa}  \int_{\Omega} \bar{g} |\nabla I|^2\, dx\,, \label{esti:gradient-inters-gw}
\end{align}
we obtain
\begin{align*}
  \frac{d}{dt} \Big[\| I_t|_{L^2}^2  + \int_{\Omega} \bar{g} |\nabla I |^2\, dx\Big] 
  \le C\,\Big(\|I_t\|_{L^2}^2 + \int_{\Omega} \bar{g} |\nabla I |^2\, dx\Big)\,.
\end{align*}
 An application of Gronwall's lemma along with \eqref{esti:gradient-inters-gw} gives: for a.e. $t\in (0,T]$
\begin{align}\label{bound_I_t_nabla_I}
\| I_t(t)\|_{L^2}^2   + \|\nabla I(t)\|_{L^2}^2 \leq  C e^{C\,t} \,.
\end{align}
 Since $I(x,t)=I_0(x)+ \displaystyle \int_{0}^{t} I_t(x,s)\,ds$, thanks to Young's inequality and \eqref{bound_I_t_nabla_I}, we have $\|I(t)\|_{L^2}^2 \le C_T \|I_0\|_{H^1}^2$ and hence
\begin{align*}
\| I\|_{L^{\infty}(0,T;H^{1})} + \| I_t\|_{L^{\infty}(0,T;L^2)} \leq C \| I_0\|_{H^1}\,.
\end{align*}
\noindent{Proof of ${\rm b)}$:} Choose $\phi \in H^1$ with $||\phi||_{H^1}\leq 1$ in \eqref{linmaina}, and use Cauchy-Schwarz inequality along with ${\rm a)}$, Lemma \ref{lem:a-priori} to obtain $\big| \langle I_{tt}, \phi  \rangle \big| 
\leq  C\,\|I_0\|_{H^1} \|\phi\|_{H^1}$ and hence 
\begin{align*}
\| I_{tt}\|_{(H^1)^\prime} \leq C \|I_0\|_{H^1}\,. 
\end{align*}
Therefore ${\rm b)}$ follows once we take square both side of the above inequality and then integrate over $(0,T)$. 
%\end{proof}
%========================================================================================================================================================
\subsection{Proof of Theorem \ref{thm:existence-uniqueness}}
In this section, we prove wellposedness of weak solution of the underlying problem via Schauder's fixed-point theorem.  To proceed further, we introduce the subspace $W_0$ of $W(0,T)$ defined by
\begin{align*}
W_0=\Big\{ & w \in W(0,T):\,  \|w\|_{L^\infty(0,T; H^1)} + \|w_t\|_{L^\infty(0,T; L^2)} \leq C\|I_0\|_{H^1}^2\,;\\
& \hspace{2cm} ~~ 0<\alpha \le w(x,t)~{\rm for ~a.e.}~(x,t)\in \Omega_T\,,~~\text{and}~~w~{\rm satisfies}~\eqref{mainc}\Big\}\,.
\end{align*}
Moreover, one can prove that $W_0$ is a non-empty, convex and weakly compact subset of $W$. Consider a mapping
\begin{align*}
\mathcal{P}:~ & W_0 \goto W_0 \\
& w\mapsto I_w\,.
\end{align*}
In order to use Schauder's fixed-point theorem on $\mathcal{P}$, we need to prove only that the mapping $\mathcal{P}:w \rightarrow I_w $ is weakly continuous from $W_0$ into $W_0$. Let
$w_k$ be a sequence that converges weakly to some $w$ in $W_0$ and let $I_k = I_{w_k}$. We have to show that $\mathcal{P}(w_k):= I_k$ converges weakly
to $\mathcal{P}(w): = I_w$.

Thanks to Lemma \ref{lem:a-priori}, one can use classical results of compact inclusion in Sobolev spaces \cite{raadams1975}, to extract subsequences $\{w_{k_n}\}$ of $\{w_k\}$ and $\{I_{k_n}\}$ of $\{I_k\}$ such that
for some $I\in W_0$, the following hold as $k\goto \infty:$
\begin{align*}
\begin{cases}
 w_{k} \longrightarrow  w \hspace{0.2cm} \text{in} \hspace{0.2cm} L^2(0,T;L^2) \hspace{0.2cm} \text{ and a.e. on } \hspace{0.2cm}  \Omega_T,\\[0.5em]
G_{\xi}\ast w_k  \longrightarrow G_{\xi}\ast w  \hspace{0.2cm} \text{in} \hspace{0.2cm}  L^2(0,T;L^2) \hspace{0.2cm}  \text{and a.e. on} \hspace{0.2cm}   \Omega_T,\\[1em]
| G_{\xi}\ast w_k |^\nu \longrightarrow | G_{\xi}\ast w |^\nu \hspace{0.2cm} \text{in} \hspace{0.2cm}  L^2(0,T;L^2) \hspace{0.2cm}  \text{and a.e. on} \hspace{0.2cm}   \Omega_T,\\[1em]
\dfrac{| G_{\xi}\ast w_k |^\nu}{ \big(M_\xi^{w_k}\big)^\nu + |G_{\xi}\ast w_k |^\nu} \rightarrow \dfrac{|G_{\xi}\ast w |^\nu}{ \big(M_\xi^{w}\big)^\nu + |G_{\xi}\ast w |^\nu} \hspace{0.2cm} {\rm in}~~L^2(0,T;L^2) \hspace{0.2cm} {\rm and~a.e.~on}~~\Omega_T\,,\\[1em]
\partial_{x_i} G_{\xi}\ast w_k  \rightarrow \partial_{x_i} G_{\xi}\ast w ~(i=1,2) \hspace{0.2cm} {\rm in}~~L^2(0,T;L^2) \hspace{0.2cm} {\rm and~a.e.~on}~~\Omega_T\,,\\[1em]
\dfrac{1}{1 + \left(\frac{|\nabla G_{\xi}\ast w_k|}{K}\right)^2} \longrightarrow \dfrac{1}{1 + \left(\frac{|\nabla G_{\xi}\ast w|}{K}\right)^2}  \hspace{0.2cm} \text{in} \hspace{0.2cm}  L^2(0,T;L^2) \hspace{0.2cm}  \text{and a.e. on} \hspace{0.2cm}   \Omega_T,\\[1.5em]
\displaystyle I_{k} \rightarrow  I\, \hspace{0.2cm} \text{weakly} *~ \text{in}~~L^{\infty}(0,T;H^1)\,,\\[1em]
\displaystyle I_{k} \rightarrow  I\, \hspace{0.2cm} \text{in}~~L^{2}(0,T; L^2)\,,\\[1em]
\partial_t I_k \rightarrow  \partial_t I\, \hspace{0.2cm} \text{weakly} * ~\text{in}~~L^{\infty}(0,T;L^2)\,,\\[1em]
\partial_{tt} I_k \rightarrow \partial_{tt}I \hspace{0.2cm} \text{weakly} *~ \text{in}~~L^{2}(0,T;(H^1)^\prime)\,.
\end{cases}
\end{align*}

The above convergence allow us to pass to the limit in the problem \eqref{linmaina} and obtain $I=\mathcal{P}(w)$.  Moreover, since the solution of \eqref{linmaina} is unique, the whole
sequence $I_k=\mathcal{P}(w_k)$ converges weakly in $W_0$ to $I=\mathcal{P}(w)$. Hence $\mathcal{P}$ is weakly continuous. Consequently, thanks to the Schauder fixed
point theorem, there exists $w \in W_0$ such that $w=\mathcal{P}(w)=I_w$.  Thus, the function $I_w$ solves the problem \eqref{maina}-\eqref{mainc}.
%==================================================================================================================
\vspace{.1cm}

\noindent{\bf Uniqueness of weak solution:}
Following the idea as in \cite{LCEvans1998}, we prove the uniqueness of weak solutions of the underlying problem \eqref{maina}-\eqref{mainc}. Let $I_{1}$ and $I_{2}$ be two weak solutions of \eqref{maina}-\eqref{mainc}.
Then, we have
\begin{align}
&I_{tt}+ I_t-\text{div} \big(g_{I_1} \nabla I\big)  = {\rm div}\big( \big(g_{I_1}-g_{I_2}\big) \nabla I_2 \big)\hspace{1.0cm}\text{in}~~\Omega_T\,, \label{eq:maina}\\
& \begin{cases} \label{eq:mainc}
 I(x,0)= 0\,,~ I_t(x,0)=0\, \hspace{3.8cm} {\rm in}~~~\Omega\,, \\
\partial_n I =0 \hspace{6.3cm} {\rm on}~~~\partial \Omega_T\,,
\end{cases}
\end{align}
where $I=I_1-I_2$. 
It suffices to show that $ I \equiv 0$ . To verify this, fix $ 0 < s < T$, and set for $i=1,2$, 
\begin{align}\label{relationvi}
v_{i}(\cdot,t)= \begin{cases}
\displaystyle \int_{t}^{s} I_{i}(\cdot, \tau)d\tau, \hspace{0.5cm} 0<t\leq s\,, \\ 
 0 \hspace{2.5cm} s \leq t < T\,.
\end{cases}
\end{align}
Note that, for $t\in (0,T)$,
\begin{align}\label{eq:fact-1}
 \begin{cases}
 \partial_t v_i(x,t)=-I_i(x,t) \quad i=1,2\,, \\
  v_{i}(\cdot,t) \in H^1\,,~~~\partial_n v_{i}=0~~~\text{on}\,\, \partial \Omega\,\,\text{in the sence of distribution}.
 \end{cases}
\end{align}
Set $v=v_1-v_2$. Then $v(\cdot,s)=0$. 
Multiplying \eqref{eq:maina} by $v$, integrating over $\Omega \times (0,s)$ along with the integration by parts formula, \eqref{eq:fact-1}, Cauchy-Schwarz inequality and the equality
\begin{align*}
 g_{I_1}\partial_t \nabla v\cdot \nabla v = \frac{1}{2} \partial_t\big(g_{I_1} |\nabla v|^2\big)-\frac{1}{2} \partial_t g_{I_1}|\nabla v|^2\,,~~{\rm and}~~
 \nabla v(x,s)=0\,,
\end{align*}
we obtain
\begin{align}\label{unique9}
&\frac{1}{2}\|I(s)\|_{L^2}^2+\int_{0}^{s}\|I(t)\|_{L^2}^2\,dt  + \frac{1}{2}\int_{\Omega} g_{I_1}(x,0) |\nabla v(x,0)|^2\, dx \nonumber  \\
&\leq \frac{1}{2} \Big|\int_{0}^{s}\int_{\Omega} |\nabla v|^2 \partial_t g_{I_1}\, dx\,dt\Big| + \int_0^s \|(g_1-g_2)(t)\|_{L^\infty} \|\nabla I_2(t)\|_{L^2}\|\nabla v(t)\|_{L^2}\,dt\,.
\end{align}
As seen in the proof of Claim \ref{claim:1}, there exist positive constants $\kappa_1, C_1>0$ such that
\begin{align*}
\kappa_1 \leq g_{I_1}\leq 1\,,\quad |\partial_t g_{I_1}|\le C_1\,.
\end{align*}
Moreover, one can use property of convolution along the fact that solution $I_i$ has positive lower bound to show that 
\begin{align*}
\|(g_{I_1}-g_{I_2})(t)\|_{L^{\infty}} \leq C(\xi, \nu, \alpha, K I_0)\|I(t)\|_{L^{2}}^\nu\,.
\end{align*}
Thus, using the above estimates in \eqref{unique9}, we have for $\nu\ge 1$
\begin{align*}
\frac{1}{2}\|I(s)\|_{L^2}^2+\int_{0}^{s}\| I(t)\|_{L^2}^2\,dt  + C \|\nabla v(0)\|_{L^2}^2 
&\le   C\, \int_0^s \big( \|\nabla v(t)\|_{L^2}^2 + \|I(t)\|_{L^2}^2\big)\,dt\,.
\end{align*}
%In the last inequality, we have used the fact that $\alpha \le 1$. 
Since $ \|v(0)\|_{L^2}^2 \le T \int_0^s \| I(t)\|_{L^2}^2\,dt $, we have 
\begin{align}\label{unique15_1}
\frac{1}{2}\|I(s)\|_{L^2}^2+\int_{0}^{s}\| I(t)\|_{L^2}^2\,dt  + C \|v(0)\|_{H^1}^2 
\le C\,\int_0^s \big( \|v(t)\|_{H^1}^2 + \|I(t)\|_{L^2}^{2} \big)\,dt\,.
\end{align}
Set
\begin{align*}
w_{i}(\cdot,t)=& \int_{0}^{t} I_{i}(\cdot,\tau)d\tau\, ; \quad w(\cdot,t)=(w_1-w_2)(\cdot,t)\,, \hspace{0.5cm} 0<t\leq T.
\end{align*}
Observe that 
\begin{align*}
 v(x,0)= w(x,s) \quad {\rm and}~~
 v(x,t)= w(x,s)-w(x,t)~~{\rm for}~~ 0<t\le s\,.
\end{align*}
Hence \eqref{unique15_1} reduces to
\begin{align}\label{unique16}
&\frac{1}{2}\|I(s)\|_{L^2}^2+\int_{0}^{s}\| I(t)\|_{L^2}^2\,dt  + C \|w(s)\|_{H^1}^2 \notag \\
& \le \tilde{C} s\,\|w(s)\|_{H^1}^2 + C\,\int_0^s \Big( \|w(t)\|_{H^1}^2 + \|I(t)\|_{L^2}^{2} \Big)\,dt\,.
\end{align}
Choose $T_1$ sufficiently small such that  $C-\tilde{C} T_1 >0$. 
Then, for $0<s\leq T_1,$ we have, from \eqref{unique16}
\begin{align}
\| I(s)\|_{L^2}^2 + \|w(s)\|_{H^1}^2 \le C \int_0^s\Big( \|w(t)\|_{H^1}^2 + \|I(t)\|_{L^2}^{2}\Big)\,dt\,. \label{unique17}
\end{align}
Consequently, an application of Gronwall's lemma then implies $ I \equiv 0$ on $[0,T_1]$.
Finally, we utilize a similar logic on the intervals $(T_1, 2T_1]$, $(2T_1,3T_1],\ldots$ step by step, and eventually deduce that
$I_{1} = I_{2}$ on $(0,T)$. This finishes the proof of Theorem \ref{thm:existence-uniqueness}.
%=====================================================================================================================
\begin{lem}
 Let $I$ be a weak solution of the problem \eqref{maina}-\eqref{mainc}, and $\beta_1:= \underset{x\in \Omega}\sup I_0(x) < \infty$.  Then 
 \begin{align}
  0<\alpha \le I(t,x)\le \beta_1 \quad {\rm for ~a.e.}~(x,t)\in \Omega_T\,. \label{boundedness-weak-solution}
 \end{align}
\end{lem}
\textbf{Proof:}
%\begin{proof}
Integrating the equation \eqref{maina} w.r.to time variable and using \eqref{mainc}, we get that
\begin{align}
 I_t + (I-I_0) -\int_0^t {\rm div}\big( g_{I}(x,s) \nabla I\big)\,ds=0 \quad \forall~(x,t)\in \Omega_T\,. \label{maind}
\end{align}
Note that, $(I-\beta)_{+} \in H^1$, where $(\cdot)_{+}$ is the truncated function defined as $(\theta)_{+}=\max\{ 0,\theta\}$. Multiplying the PDE \eqref{maind} by $(I-\beta_1)_{+}$
and then integrating over $\Omega$ to have
\begin{align*}
 \frac{1}{2} \frac{d}{dt} \int_{\Omega} |(I-\beta_1)_{+}|^2\,dx + \int_{\Omega} (I-I_0)(I-\beta_1)_{+}\,dx + \int_0^t \int_{\{I\ge \beta_1\}} g_{I}(x,s)|\nabla I|^2\,dx\,ds=0\,.
\end{align*}
Observe that, $g_{I} \ge 0$ and $(I-I_0)(I-\beta_1)_{+} \ge 0$. Thus, we have $\frac{d}{dt} \int_{\Omega} |(I-\beta_1)_{+}|^2\,dx\le 0$. Again, since $I_0 \le \beta_1$, we obtain 
$\int_{\Omega} |(I-\beta_1)_{+}|^2\,dx \le 0$ for a.e. $t\in [0,T]$. Therefore,  $I(x,t)\le \beta_1$ for a.e. $(x,t)\in \Omega_T$. 
\vspace{.1cm}

Similarly, multiplying the equation \eqref{maind} with $(I-\alpha)_{-}\in H^1$ and then integrating over $\Omega$ to conclude that $0<\alpha \le I(x,t)$ for a.e. $(x,t)\in \Omega_T$, where 
$(\cdot)_{-}$ is the truncated function defined as $(\theta)_{-}=\min\{ 0,\theta\}$. Hence \eqref{boundedness-weak-solution} holds true. This completes the proof. 
 
%\end{proof}
%===================================================================================================================================
\section{Numerical Implementation}
\label{sec:numerical}
To solve the present model numerically, we choose an explicit finite difference scheme, which is the most straightforward option for solving a hyperbolic PDE.\\
(a). Discretize the time domain using a step $\tau$ and the space domain using a step $ h.$ Denote $I^n_{i,j}=I(x_i,y_j,t_n)$ where $x_i=ih, \hspace{0.2cm} i=0,1,2...,N;$
$y_j=jh, \hspace{0.2cm} j=0,1,2...,M;$ $t_n=n\tau,\hspace{0.2cm}  n=0,1,2...$  where $n$ is the number of iterations and $M \times N$ is the size of the image.\\
(b). Boundary conditions are given as:
$I_{-1,j}^n=I_{0,j}^n,I_{N+1,j}^n=I_{N,j}^n, \hspace{0.2cm} I_{i,-1}^n=I_{i,0}^n,I_{i,M+1}^n=I_{i,M}^n.$	\\	
(c). The approximation of derivative terms are given as follows:
\begin{center}
$
\begin{array}{lll}
\displaystyle\frac{\partial I_{i,j}^n}{\partial t} & \approx & \displaystyle\frac{I_{i,j}^{n+1}-I_{i,j}^n}{\tau},
\displaystyle\frac{\partial^2 I_{i,j}^n}{\partial t^2} \approx \displaystyle\frac{I_{i,j}^{n+1}-2I_{i,j}^n+I_{i,j}^{n-1}}{\tau ^2},\\[0.8cm]
\nabla_x I_{i,j}^n &\approx& \displaystyle\frac{I_{i+h,j}^n-I_{i-h,j}^n}{2h},
\nabla_y I_{i,j}^n \approx \displaystyle\frac{I_{i,j+h}^n-I_{i,j-h}^n}{2h},\\[0.8cm]
|\nabla I_{i,j}^n| & \approx & \sqrt{(\nabla_x I_{i,j}^n)^2 + (\nabla_y I_{i,j}^n)^2}.
\end{array}
$
\end{center}
(d). The discrete form of the proposed model \eqref{maina} could be written as follows:
\begin{align*}
&(1+\gamma \tau)I_{i,j}^{n+1}=(2+\gamma \tau)I_{i,j}^n -I_{i,j}^{n-1} + {\tau ^2} \big\{ \nabla_x \left( g_{i,j}^n \nabla_x I_{i,j}^n  \right) 
 + \nabla_y \left( g_{i,j}^n \nabla_y I_{i,j}^n \right) \big\} ,
\end{align*}
where
\begin{align*}
g_{i,j}^n=b(s^{n}_{i,j})\cdot  \dfrac{1}{1+ \left(\frac{|\nabla G_{\xi} \ast I^n_{i,j} |}{K} \right)^2 },
\end{align*}
with the conditions,
\begin{align*}
\begin{split}
I_{i,j}^0 &=I_0(ih,jh),   \hspace{1.5cm}	0 \leq i \leq N, 0 \leq j \leq M,
\\
I_{i,j}^1 &=I_{i,j}^0,   \hspace{2.5cm}	0 \leq i \leq N, 0 \leq j \leq M.
\end{split}
\end{align*}

Apart from the discretization of \eqref{maina}-\eqref{mainc}, we need to specify a stopping criterion for the convergence of the numerical simulation process.  For this, we start the simulation with the initial value $I_0$ and utilize the system \eqref{maina} repeatedly, resulting in a family of smoother images ${I(x,t)};t>0$, which represents filtered versions of $I_0$.  And then we stop the noise elimination process after getting the best PSNR value of the restored image.

%\iffalse
\section{Experiment Results and Discussion}
\label{sec:Results}
This section displayed the performance of the present model in terms of visual quality and quantitative results. We compare the despeckling result of the proposed model using three standard test images corrupted by the multiplicative speckle noise with a different number of looks~($L$). We have artificially added multiplicative speckle 
noise level $ L = \{ 1,3,5,10,33 \} $ by using our MATLAB program. All the numerical tests are performed under windows $7$ and MATLAB version $R2018b$ running on a desktop with an Intel Core $i5$ dual-core CPU at $2.53$ GHz with $4$ GB of memory. Image denoising using the present model has been compared with the Shan model \cite{shan2019multiplicative}.
%Since the present model is inspired by the Shan model as well as the model is a very latest diffusion-based model, it makes sense for the comparison. 
%Apart from the artificially noisy images, we test our model on some real SAR images corrupted with different level of speckle noise.
In this process, the considered existing model is discretized using the same explicit numerical scheme as in the proposed model. We choose an uniform time step size $\tau = 0.2$ and $\xi=1$ for each models. Details of the other parameter values for the numerical computation are given in the right-hand of Table \ref{tab:psnr_ssim_parameter}.

\subsection{Image quality measurement}
\label{sec:quality}
Since the proposed model is claimed to be an improvement over the existing diffusion models, our main aim is to compare the edge detection and denoising results, in terms of both qualitative and quantitative
measures.  For each experiment, we compute the values of the two standard parameters peak signal to noise ratio (PSNR)\cite{gonzalez2002digital} and Structural similarity index (SSIM)\cite{wang2004image} for the quantitative comparison with the other existing model. A higher numerical value of PSNR and SSIM suggests that the reconstructed image is closer to the noise-free image. The considered parameters are defined as follows: \\ 
(a). PSNR can measure the match between the clean and denoised data,
\begin{align*}
\text{PSNR} = 10\, \text{log}_{10} \left(\frac{\text{max}(I)^2}{\frac{1}{\text{MN}} \sum\limits_{i=1}^\text{M} \sum\limits_{j=1}^\text{N} (I(i,j)-I_t(i,j))^2  }\right).
\end{align*}
Here $I$ denotes the clean image of size $\text{M}\times \text{N}$ and max($I$) is the maximum possible pixel value of $I,$ and $I_t$ denotes the denoised image at a certain time $t.$ \newline
(b). SSIM is used to calculate the similarity between structure of clean and reconstructed image and can be given as,
\begin{align*}
\text{SSIM}(x,y)=\frac{(2\mu_x \mu_y +k_1)(2\sigma_{xy} +k_2)}{({\mu_x}^2 +{\mu_y}^2 +k_1)({\sigma_x}^2 + {\sigma_y}^2 +k_2)}.
\end{align*}
Here $\mu_x, \mu_y, {\sigma_x}^2, {\sigma_y}^2, \sigma_{xy} $ are the average, variance and covariance of $x$ and $y$, respectively. $k_1$ and $k_2$ are the variables to stabilize the division with weak denominator.\\
(c). Other typical qualitative measures have also been computed in terms of the ratio image, which can be defined as the point-by-point ratio between the degraded and the despeckled image \cite{argenti2013tutorial}. Apart from the ratio image, we also compute the 2D contour plot, 3D surface plot for the better visualization of the computational result for the proposed model as well as for the other discussed models.
\iffalse
%=========================================================================================================================
\begin{figure}[]
       \centering
       \begin{subfigure}[b]{0.21\textwidth}           
           \includegraphics[scale=0.18]{boat}           
                \caption{Boat}
                \label{fig:boat_clear}
       \end{subfigure}%
              \begin{subfigure}[b]{0.21\textwidth}           
                \includegraphics[scale=0.36]{brick}               
                \caption{Brick}
                \label{fig:brick_clear}
       \end{subfigure}% 
      \begin{subfigure}[b]{0.21\textwidth}           
                \includegraphics[scale=0.309]{circle}               
                \caption{Circle}
                \label{fig:circle_clear}
       \end{subfigure}%
       
 \caption{Test Images: (a) Natural Image, (b) Texture Image, (c) Synthetic Image.}\label{fig:all_clear_images}
\end{figure}
\fi
%============================================================================================================================
\subsection{Computational Results \& Discussion}
\label{sec:natural}
\begin{figure}
       \centering
       
       \begin{subfigure}[b]{0.24\textwidth}           
           \includegraphics[scale=0.2]{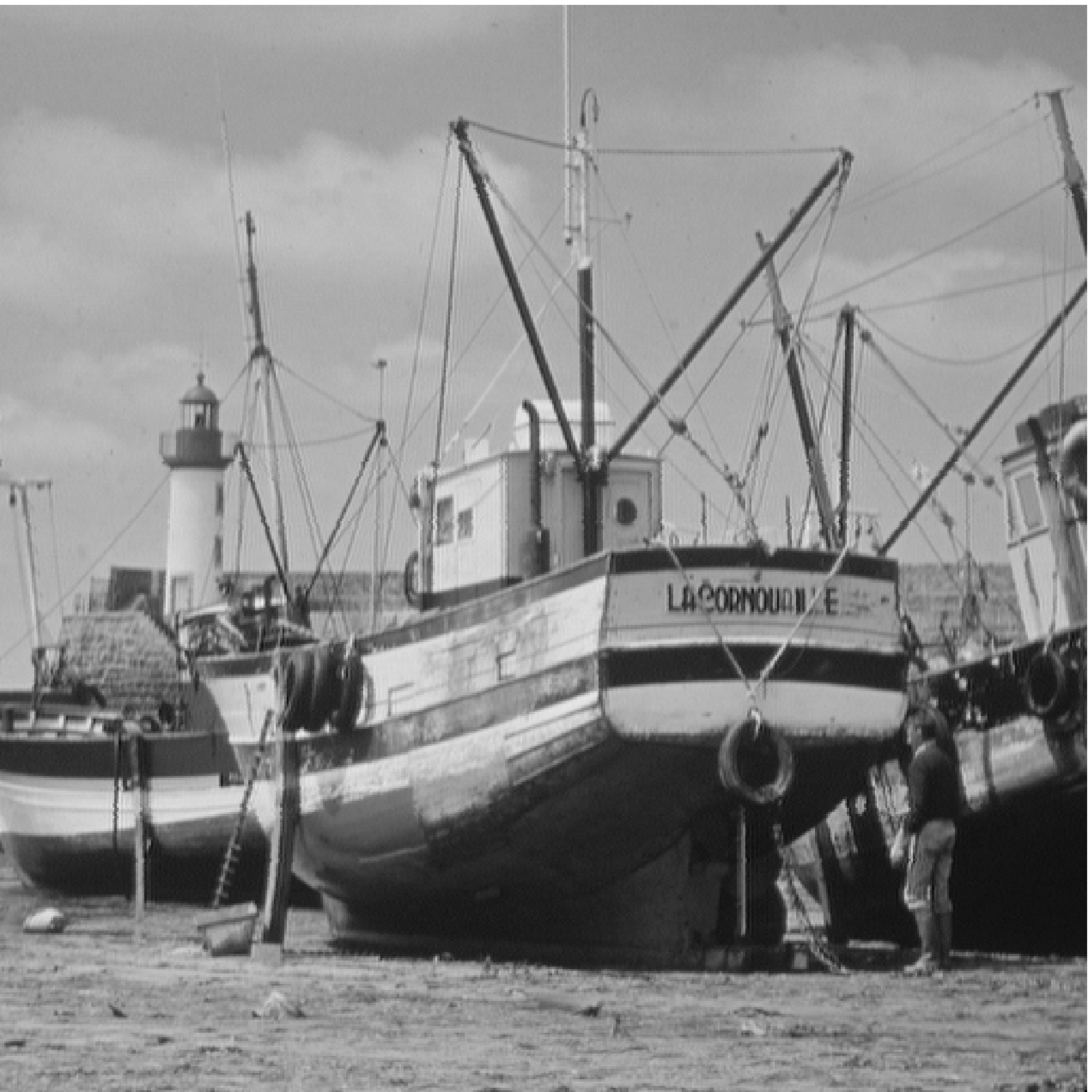}           
                \caption{Original}
                \label{fig1:boat}
       \end{subfigure}
       \begin{subfigure}[b]{0.24\textwidth}           
                \includegraphics[scale=0.2]{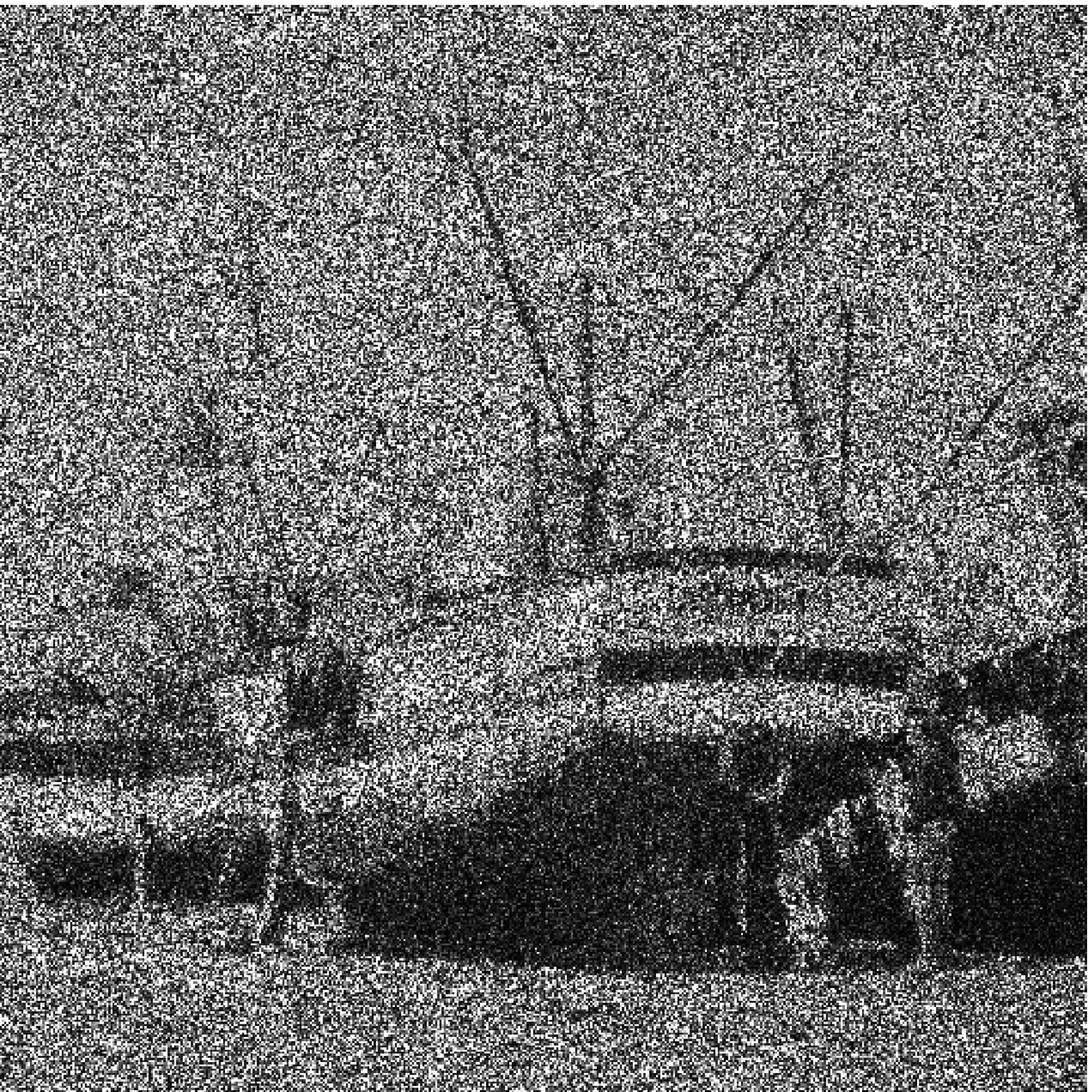}               
                \caption{Noisy}
                \label{fig:boat_1}
       \end{subfigure}% 
      \begin{subfigure}[b]{0.24\textwidth}           
                \includegraphics[scale=0.2]{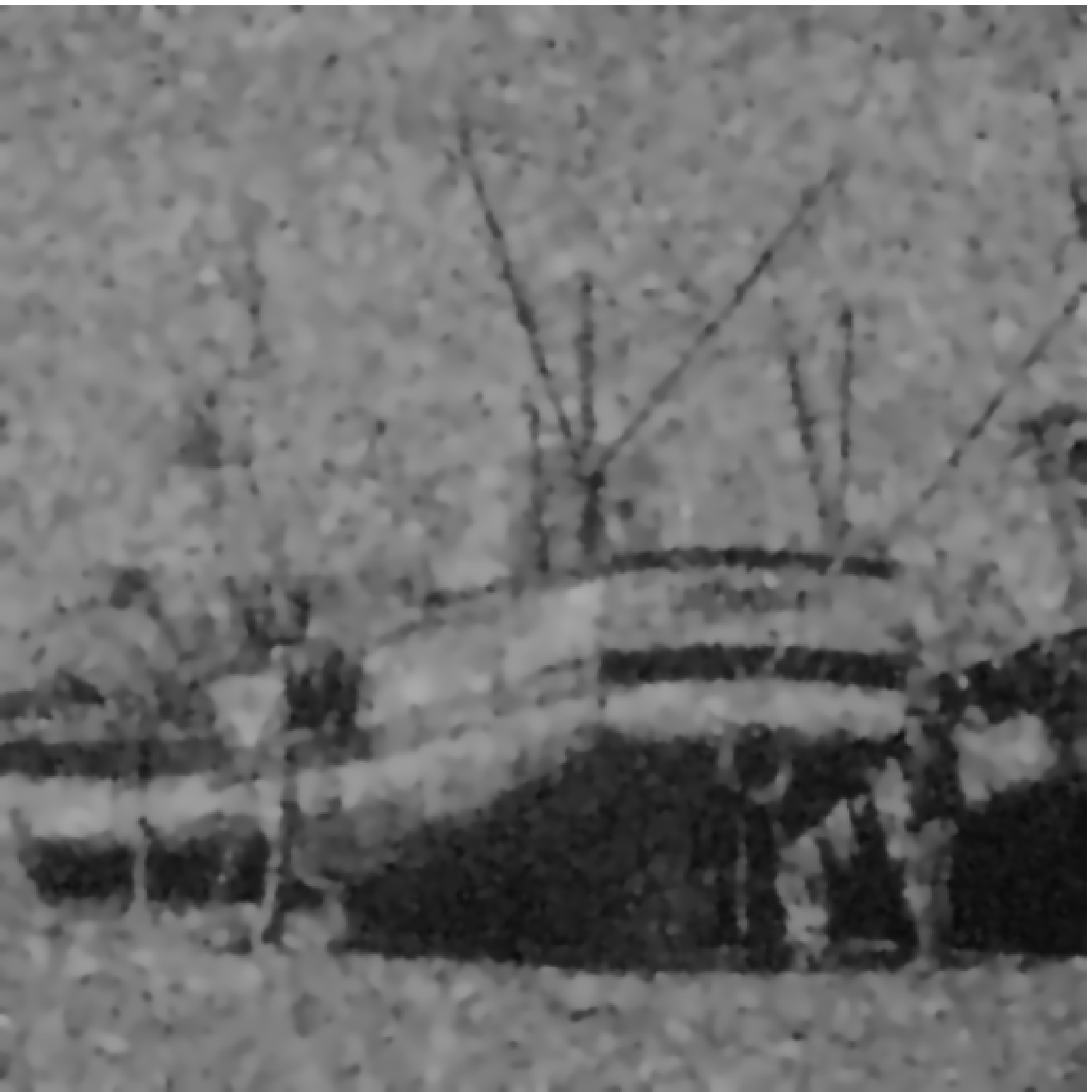}               
                \caption{Shan}
                \label{fig:boat_1_zzdb}
       \end{subfigure}% 
        \begin{subfigure}[b]{0.24\textwidth}           
                \includegraphics[scale=0.2]{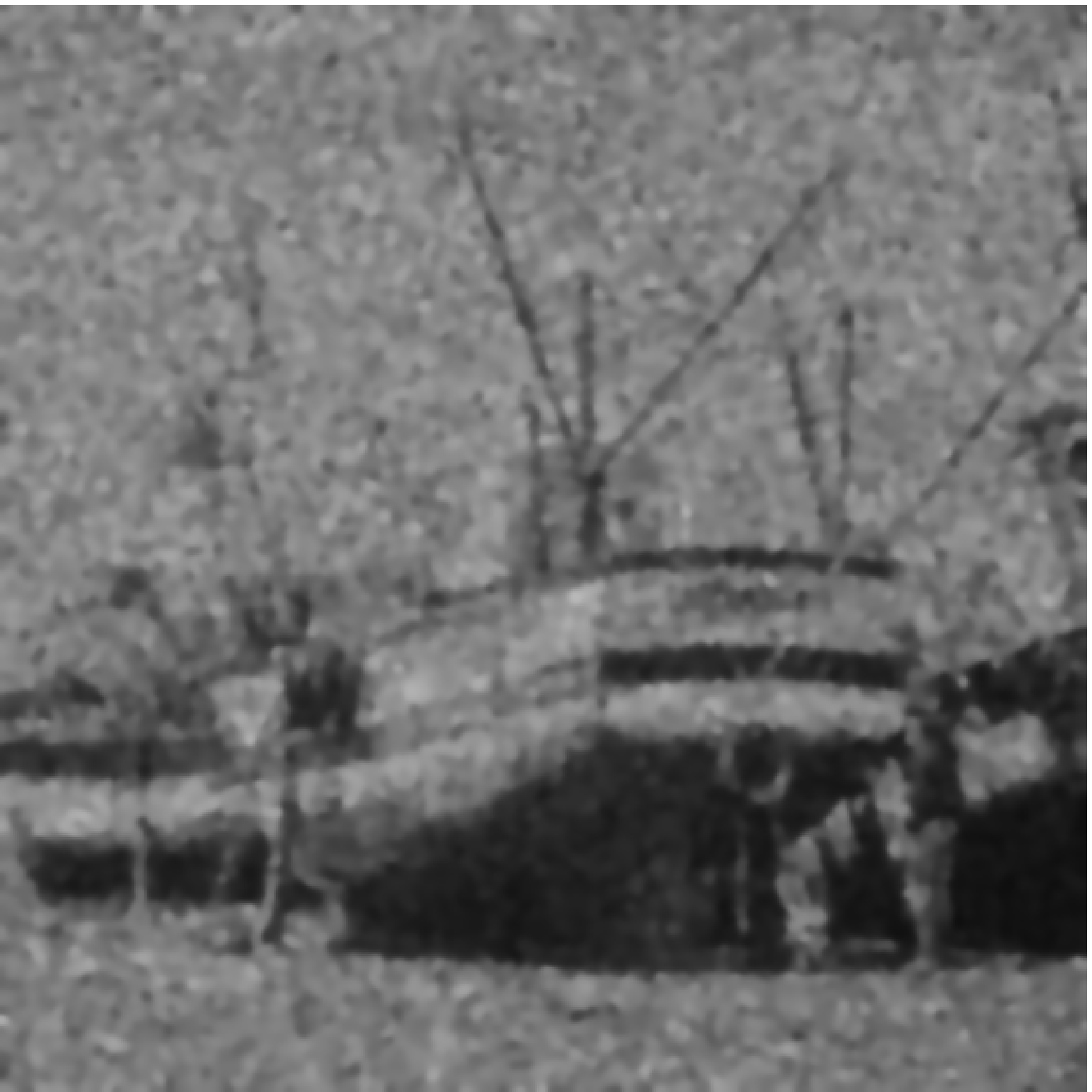}               
                \caption{Proposed}
                \label{fig:boat_1_tdm}
       \end{subfigure}
     
\caption{Image corrupted with speckle look L=1 and restored by different models.}\label{boat_1}
\end{figure}

In figure \ref{boat_1}, we represent the restored results of a Boat image (Natural Image) which is contaminated by multiplicative speckle noise with  $L=1$.  From the visual quality of the restored images, it is easy to perceive that Shan model leaves some spikes in the restored images but the results computed by the present model is more apparent than the results of Shan model.

\begin{figure}
       \centering
     
       \begin{subfigure}[b]{0.24\textwidth}           
           \includegraphics[scale=0.4]{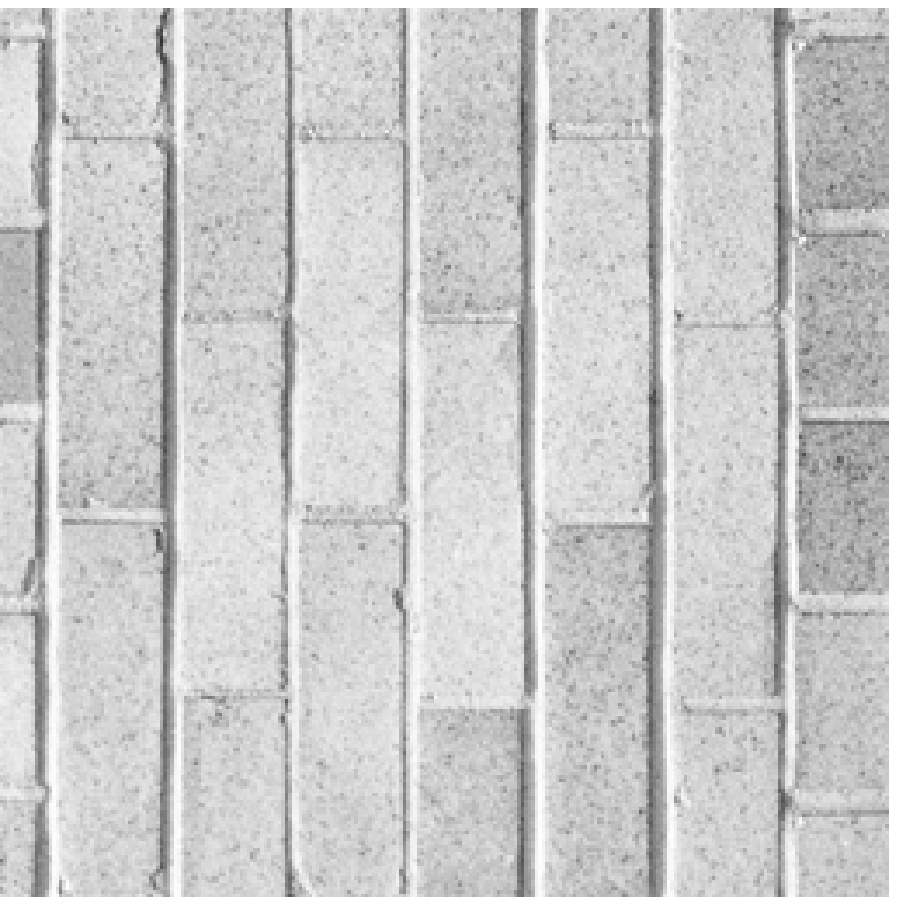}           
                \caption{Original}
                \label{fig1:brick}
       \end{subfigure}
       \begin{subfigure}[b]{0.24\textwidth}           
                \includegraphics[scale=0.4]{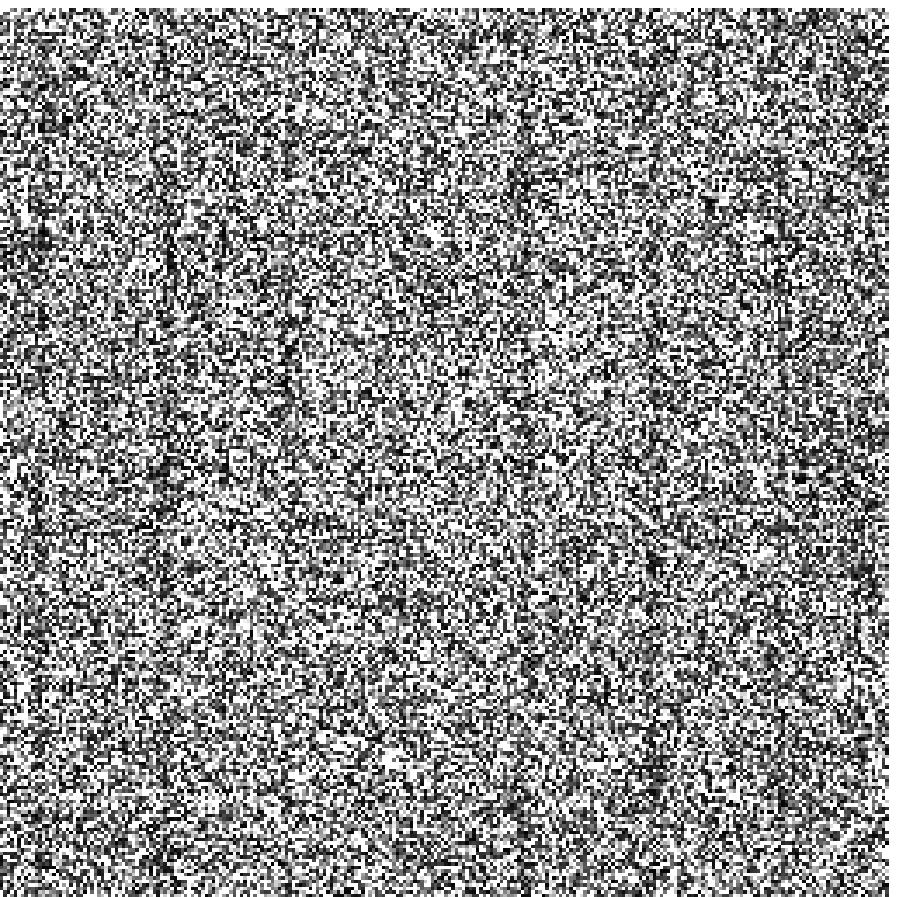}               
                \caption{Noisy}
                \label{fig:brick_1}
       \end{subfigure}% 
       \begin{subfigure}[b]{0.24\textwidth}           
                \includegraphics[scale=0.4]{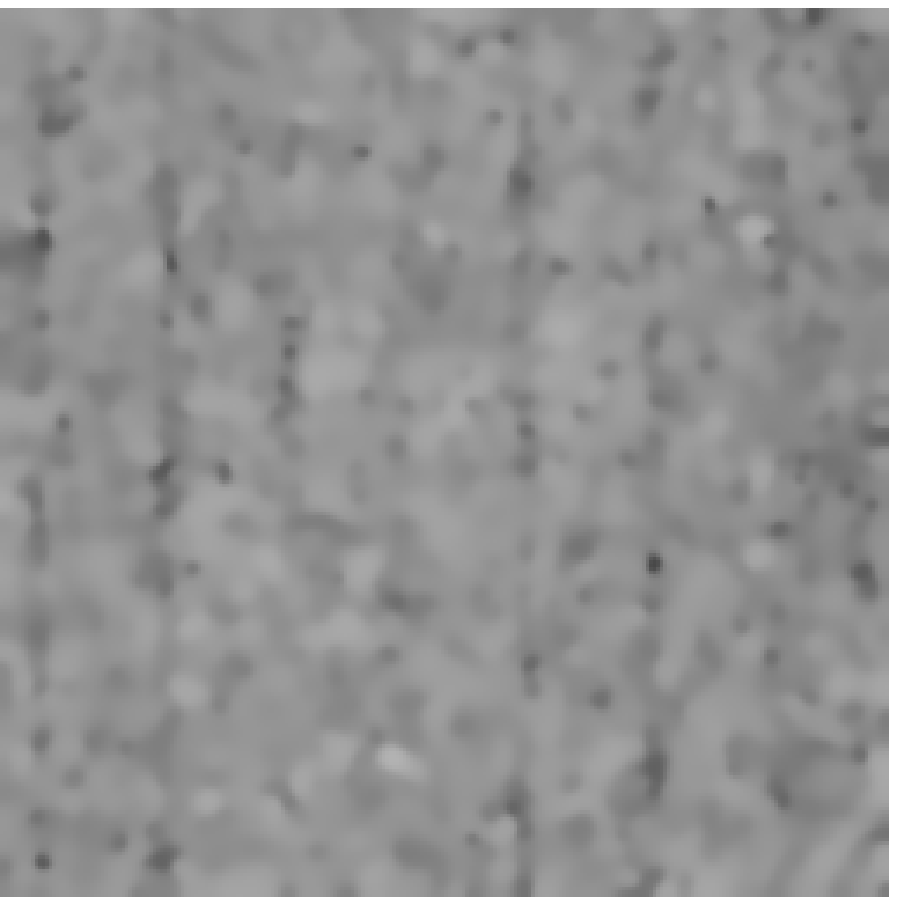}               
                \caption{Shan}
                \label{fig:brick_1_zzdb}
       \end{subfigure}% 
        \begin{subfigure}[b]{0.24\textwidth}           
                \includegraphics[scale=0.4]{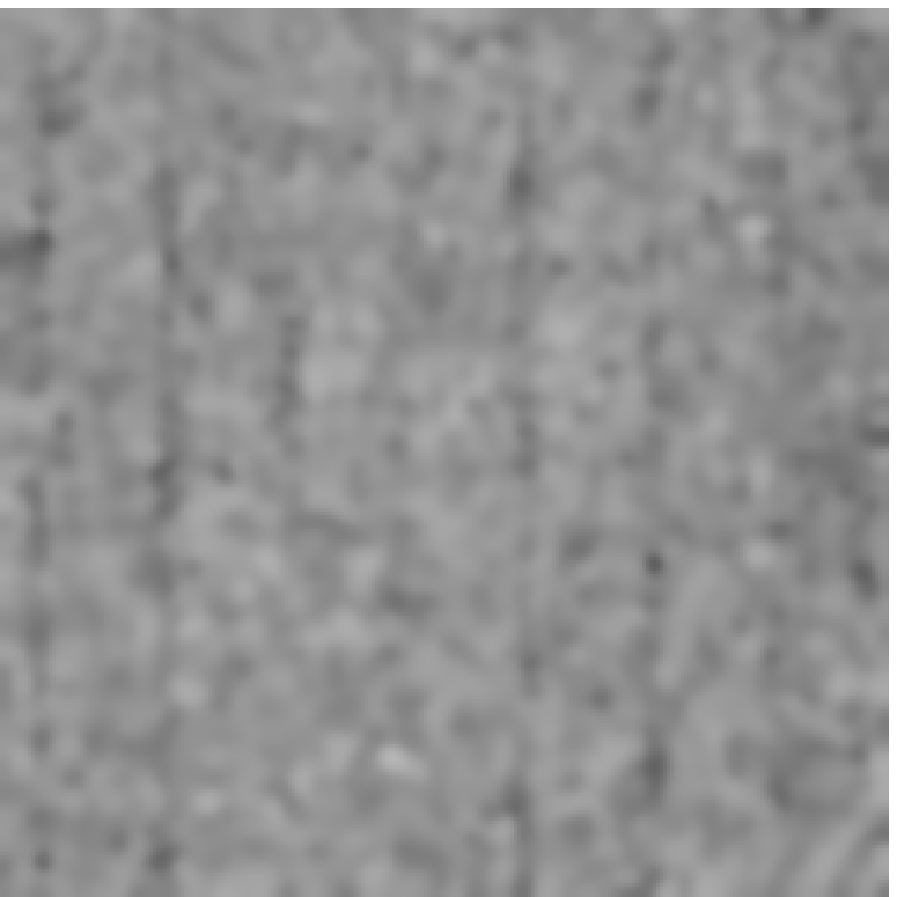}               
                \caption{Proposed}
                \label{fig:brick_1_tdm}
       \end{subfigure}
     
\caption{Image corrupted with speckle look L=1 and restored by different models.}\label{brick_1}
\end{figure}

\begin{figure}
       \centering
    
       \begin{subfigure}[b]{0.24\textwidth}           
           \includegraphics[scale=0.4]{brick}           
                \caption{Original}
                \label{fig2:brick}
       \end{subfigure}
       \begin{subfigure}[b]{0.24\textwidth}           
                \includegraphics[scale=0.4]{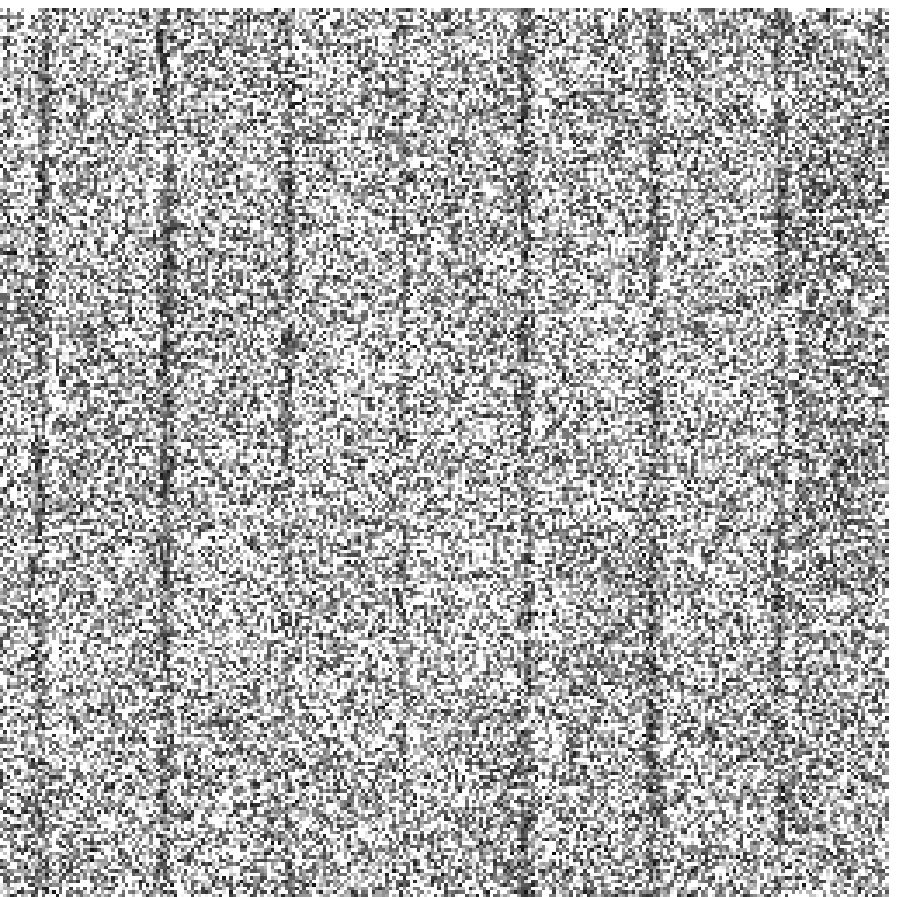}               
                \caption{Noisy}
                \label{fig:brick_3}
       \end{subfigure}% 
         \begin{subfigure}[b]{0.24\textwidth}           
                \includegraphics[scale=0.4]{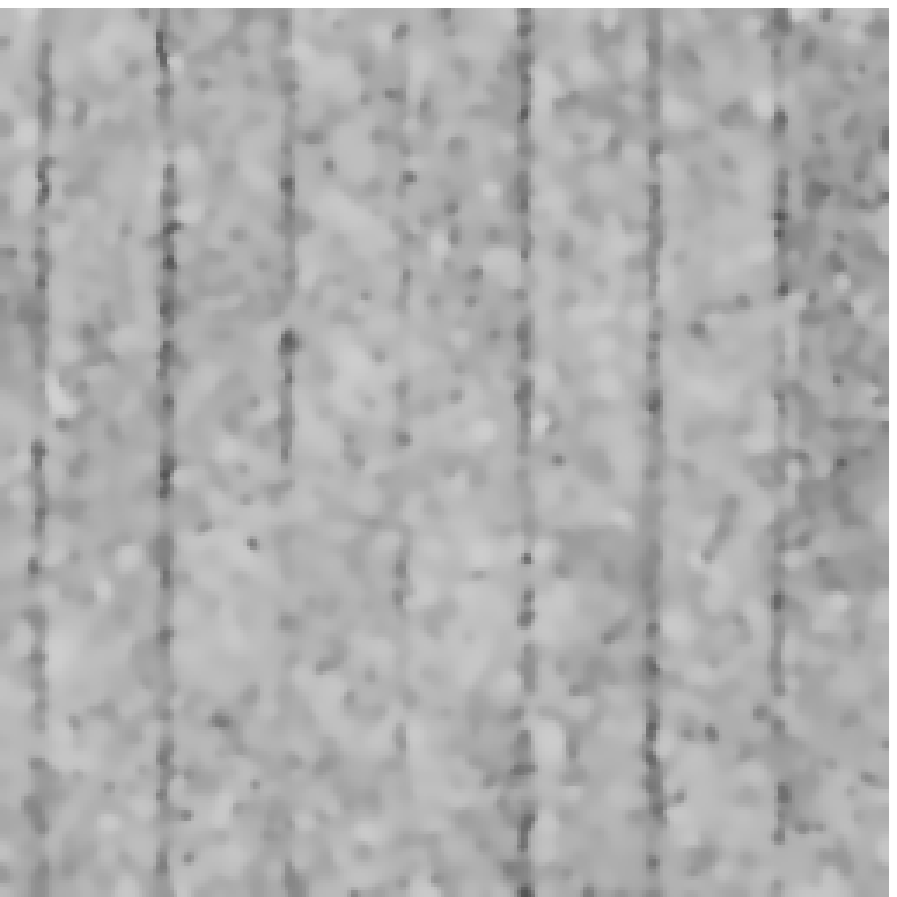}               
                \caption{Shan}
                \label{fig:brick_3_zzdb}
       \end{subfigure}% 
        \begin{subfigure}[b]{0.24\textwidth}           
                \includegraphics[scale=0.4]{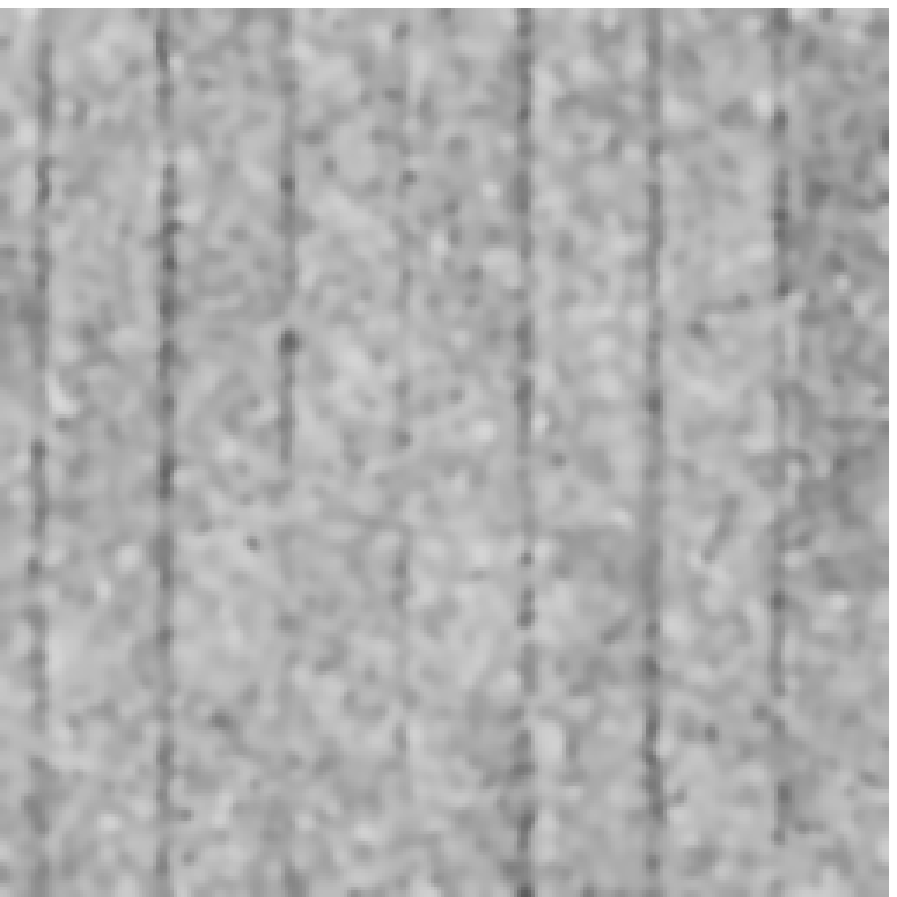}               
                \caption{Proposed}
                \label{fig:brick_3_tdm}
       \end{subfigure}
     
\caption{Image corrupted with speckle look L=3 and restored by different models.}\label{brick_3}
\end{figure}

In figures \ref{brick_1}-\ref{brick_3}, we describe the reconstructed results of a Brick image (Texture Image) which is corrupted by speckle noise with  $L=\lbrace1,3 \rbrace$.  From the figures, it is easy to see that the result computed by the present model is more apparent as well as less blurry than the Shan model.  

To check the more reconstruction capability of the present model in figures \ref{circle_1}- \ref{fig:circle_10_3d} illustrate the qualitative results of a Circle image (Synthetic Image) which is corrupted by  speckle noise with  $L = \lbrace 1,3,5,10 \rbrace$. In the images \ref{circle_1}- \ref{circle_5} we demonstrate the despeckling images by the present model and Shan model when the image is corrupted by the noise level $L=\lbrace1,3,5\rbrace$. From these figures, we easily visualize the performance of the present model.

In figure \ref{fig:circle_10_ratio}, we represent the restored image along with their ratio image for a better comparison of the qualitative result. From the figures \ref{fig:circle_10_zzdb}-\ref{fig:circle_10_tdm} it can be easily concluded that the present model gives promising results in terms of image despeckling than the Shan model. Figure \ref{fig:circle_ratio} represent the ratio image for the clear circle image \ref{fig4:circle}.Figure \ref{fig:circle_ratio} indicates that ratio image for the clear image has no background information. From the figures \ref{fig:circle_10_zzdb_ratio}- \ref{fig:circle_10_tdm_ratio} we can see that ratio image corresponding to the present model has very less background information. Which confirms that the present model works better in terms of edge preservation than the Shan model.

To more visualize the noise removal ability in figures \ref{fig:circle_10_cont}-\ref{fig:circle_10_3d} we illustrate the contour maps and 3D surface plots corresponding to the images \ref{fig4:circle}-\ref{fig:circle_10_tdm}. One can easily observe that from the contour maps, and 3D surface plots, Shan model left some speckles in the homogeneous regions, but the present model produces fewer artifacts with better edge preservation.

Along with the qualitative comparison, the quantitative results in terms of PSNR and SSIM values are displayed in Table \ref{tab:psnr_ssim_parameter}.  The highest values of PSNR and SSIM for each noise level clearly shows that the suggested model is better than the Shan model.
\iffalse
Apart from the results for artificially noisy images, in the figure \ref{fig:real_images} we display the restored results for real SAR images. Apply the present model on four different real images : $(i)$ \ref{fig:real2_noisy}: Space Radar Image of Kilauea, Hawaii \cite{dataset_kilauea}, $(ii)$ \ref{fig:real3_noisy}:  SAR image of KOMPSAT/Arirang-5 of a portion of the Himalayan Arc \cite{dataset_eoportal}, $(iii)$ \ref{fig:real4_noisy}:  High-resolution SAR image of Prague, Czech Republic \cite{dataset_eoportal}, $(iv)$ \ref{fig:real1_noisy}: One look radar image \cite{dataset_european}. 
\fi
%\iffalse
%======================================================================================================================
%\iffalse

%======================================================================================================
%======================================================================================================
\iffalse
\begin{figure}
       \centering
    
       \begin{subfigure}[b]{0.24\textwidth}           
           \includegraphics[scale=0.4]{brick}           
                \caption{Original}
                \label{fig:brick}
       \end{subfigure}
       %
       \begin{subfigure}[b]{0.24\textwidth}           
                \includegraphics[scale=0.4]{brick_speckle_look_5}               
                \caption{Noisy}
                \label{fig:brick_5}
       \end{subfigure}% 
            %
         \begin{subfigure}[b]{0.24\textwidth}           
                \includegraphics[scale=0.4]{brick_speckle_look_5_shan}               
                \caption{Shan}
                \label{fig:brick_5_zzdb}
       \end{subfigure}% 
        \begin{subfigure}[b]{0.24\textwidth}           
                \includegraphics[scale=0.4]{brick_speckle_look_5_tdm_rev}               
                \caption{Proposed}
                \label{fig:brick_5_tdm}
       \end{subfigure}
 \caption{Image corrupted with speckle look L=5 and restored by different models.}\label{brick_5}
\end{figure}
\fi
%======================================================================================================
\begin{figure}
       \centering
       
          \begin{subfigure}[b]{0.25\textwidth}           
                \includegraphics[scale=0.37]{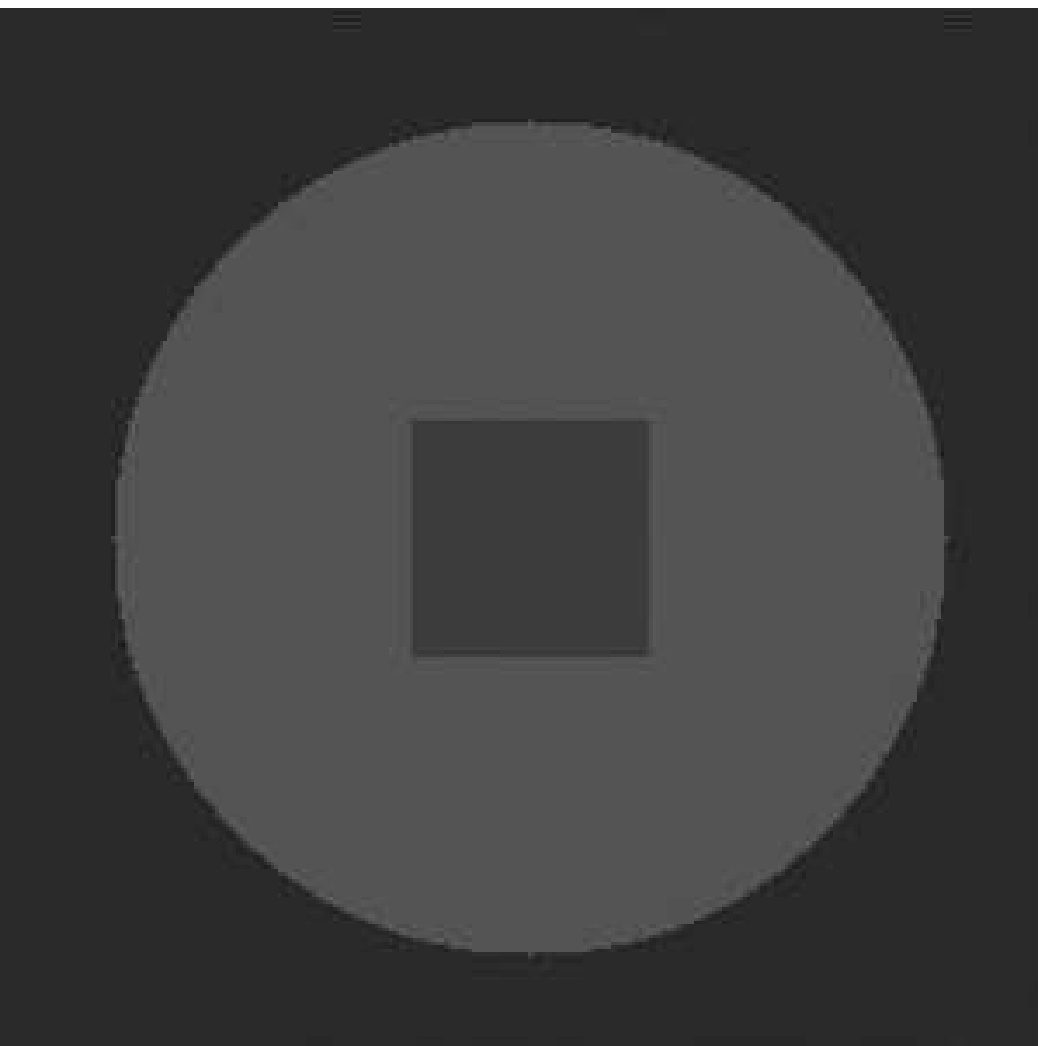}               
                \caption{Original}
                \label{fig1:circle}
       \end{subfigure}%
       \begin{subfigure}[b]{0.25\textwidth}           
                \includegraphics[scale=0.37]{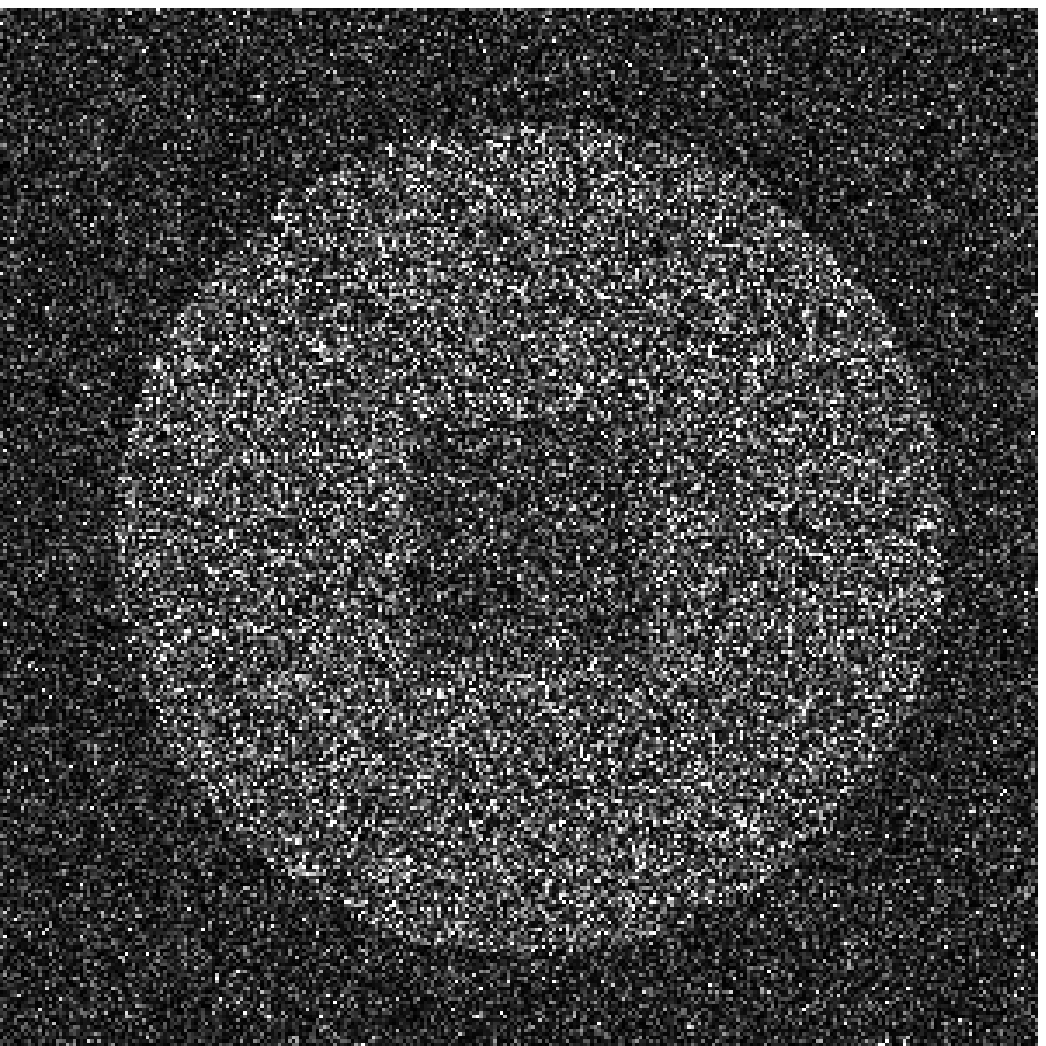}               
                \caption{Noisy}
                \label{fig:circle_1}
       \end{subfigure}% 
         \begin{subfigure}[b]{0.25\textwidth}           
                \includegraphics[scale=0.37]{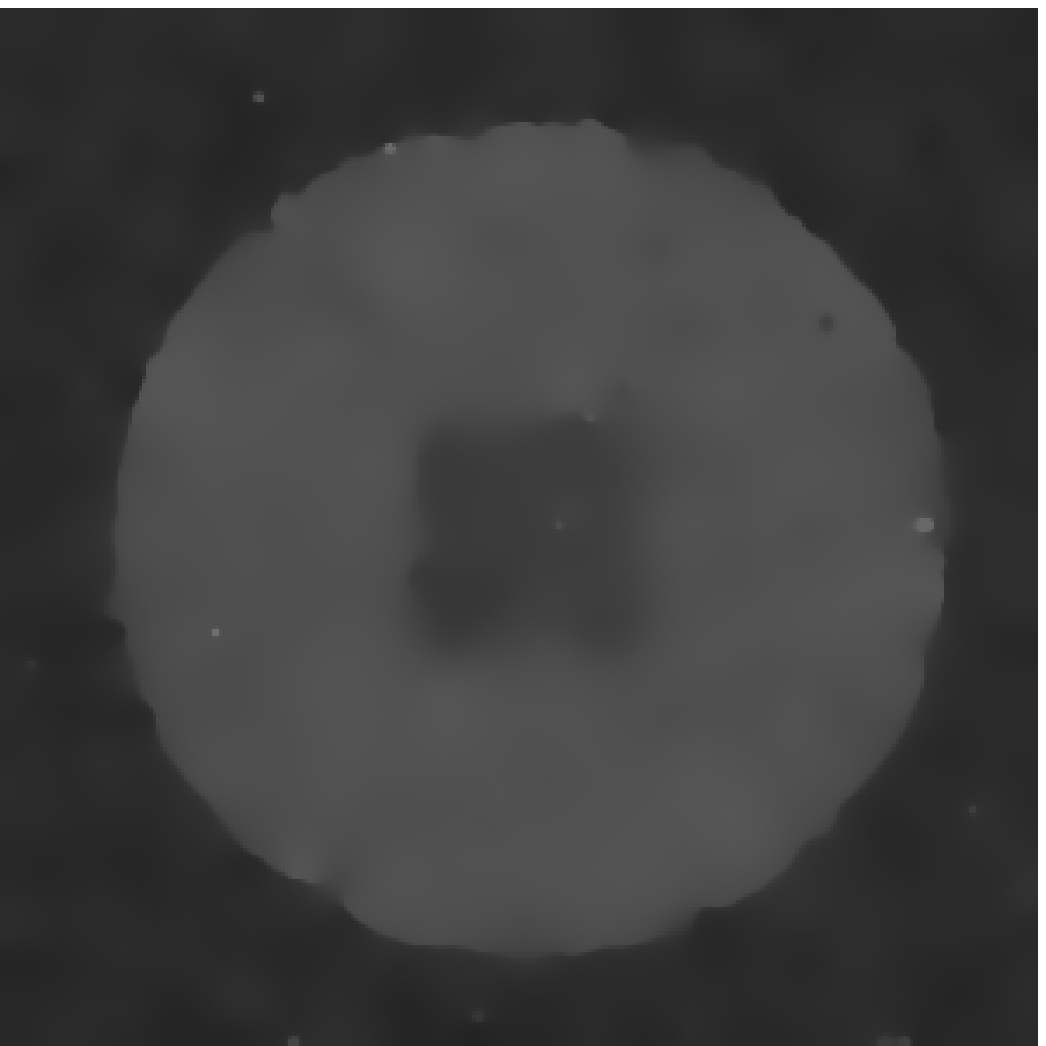}               
                \caption{Shan}
                \label{fig:circle_1_zzdb}
       \end{subfigure}% 
        \begin{subfigure}[b]{0.25\textwidth}           
                \includegraphics[scale=0.37]{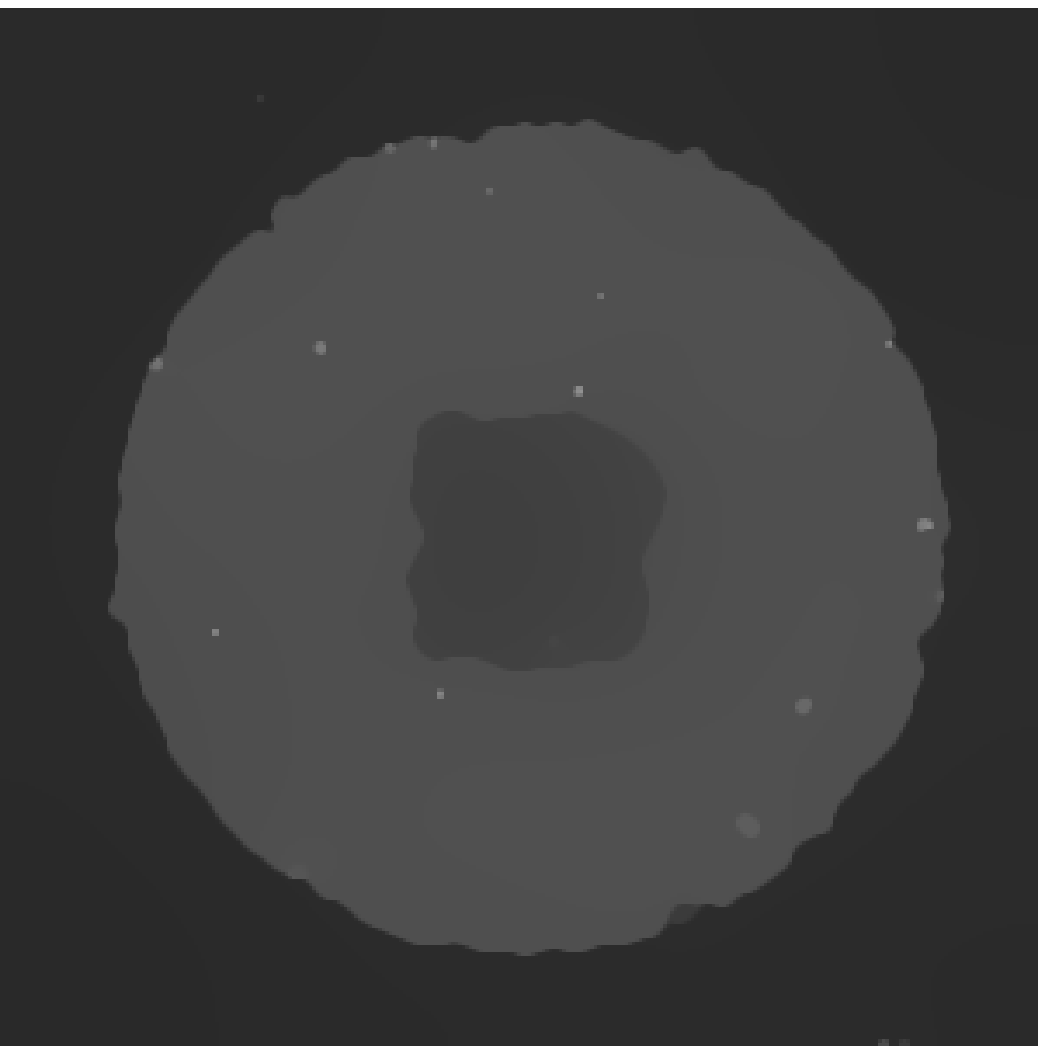}               
                \caption{Proposed}
                \label{fig:circle_1_tdm}
       \end{subfigure}
     
\caption{Image corrupted with speckle look L=1 and restored by different models. }\label{circle_1}
\end{figure}
%======================================================================================================
%======================================================================================================
\begin{figure}
       \centering
              \begin{subfigure}[b]{0.25\textwidth}           
                \includegraphics[scale=0.37]{circle}               
                \caption{Original}
                \label{fig2:circle}
       \end{subfigure}%
       \begin{subfigure}[b]{0.25\textwidth}           
                \includegraphics[scale=0.37]{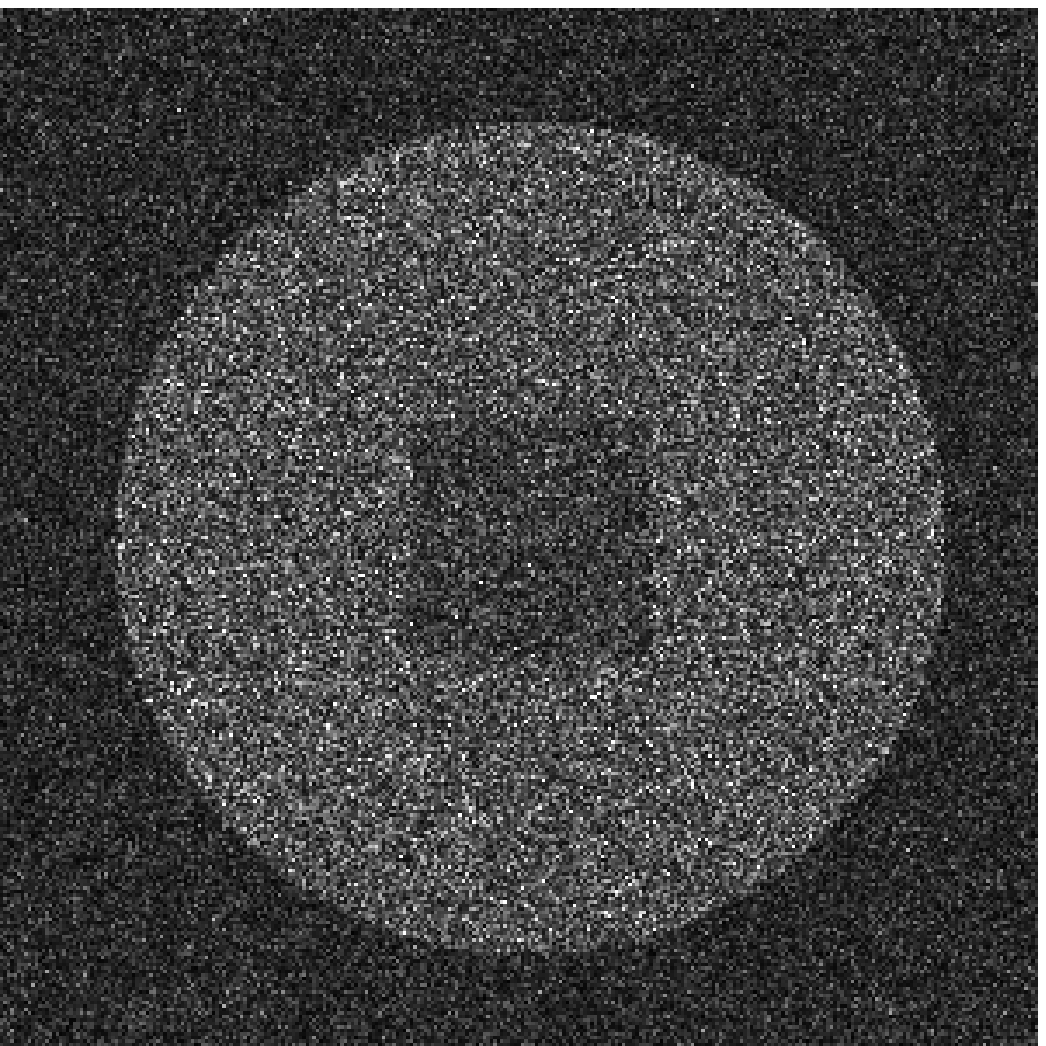}               
                \caption{Noisy}
                \label{fig:circle_3}
       \end{subfigure}% 
         \begin{subfigure}[b]{0.25\textwidth}           
                \includegraphics[scale=0.37]{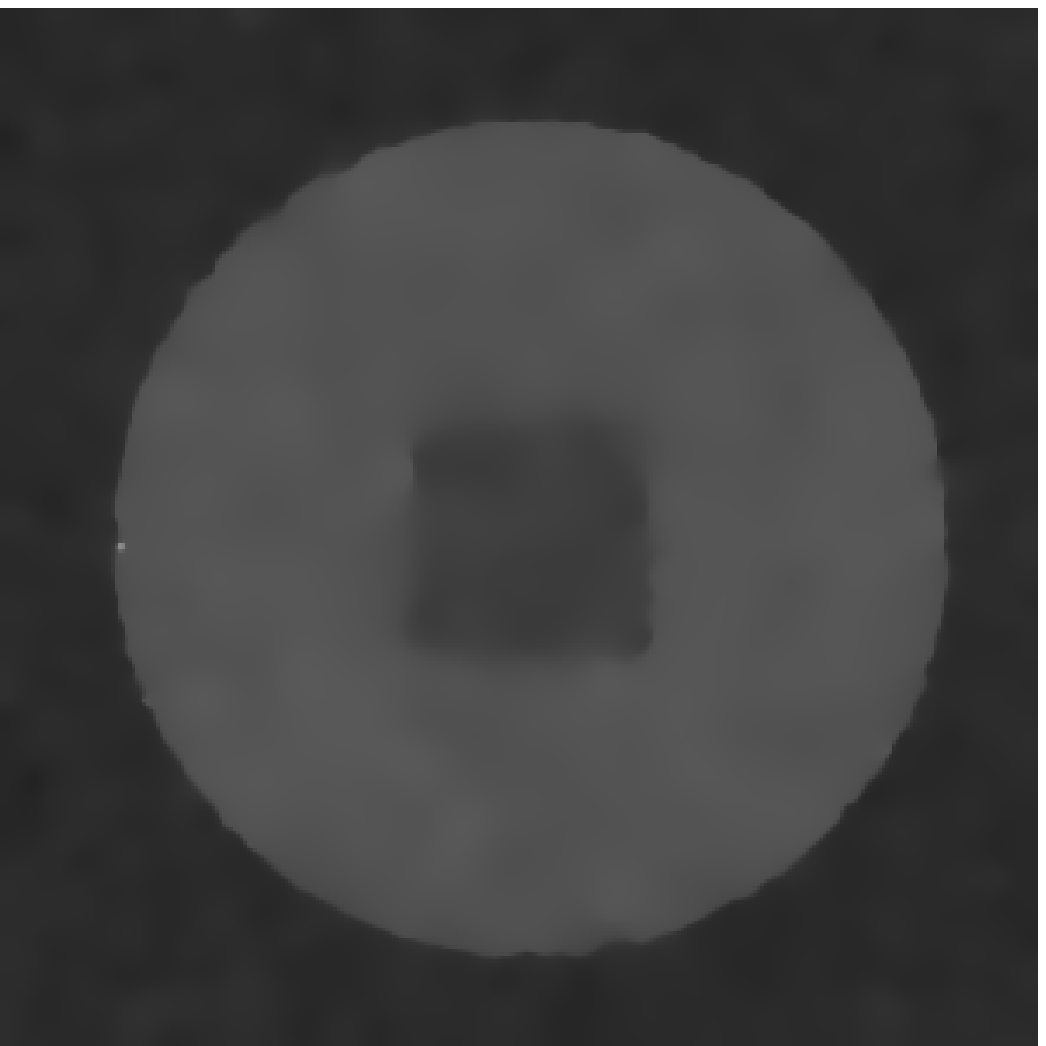}               
                \caption{Shan}
                \label{fig:circle_3_zzdb}
       \end{subfigure}% 
        \begin{subfigure}[b]{0.25\textwidth}           
                \includegraphics[scale=0.37]{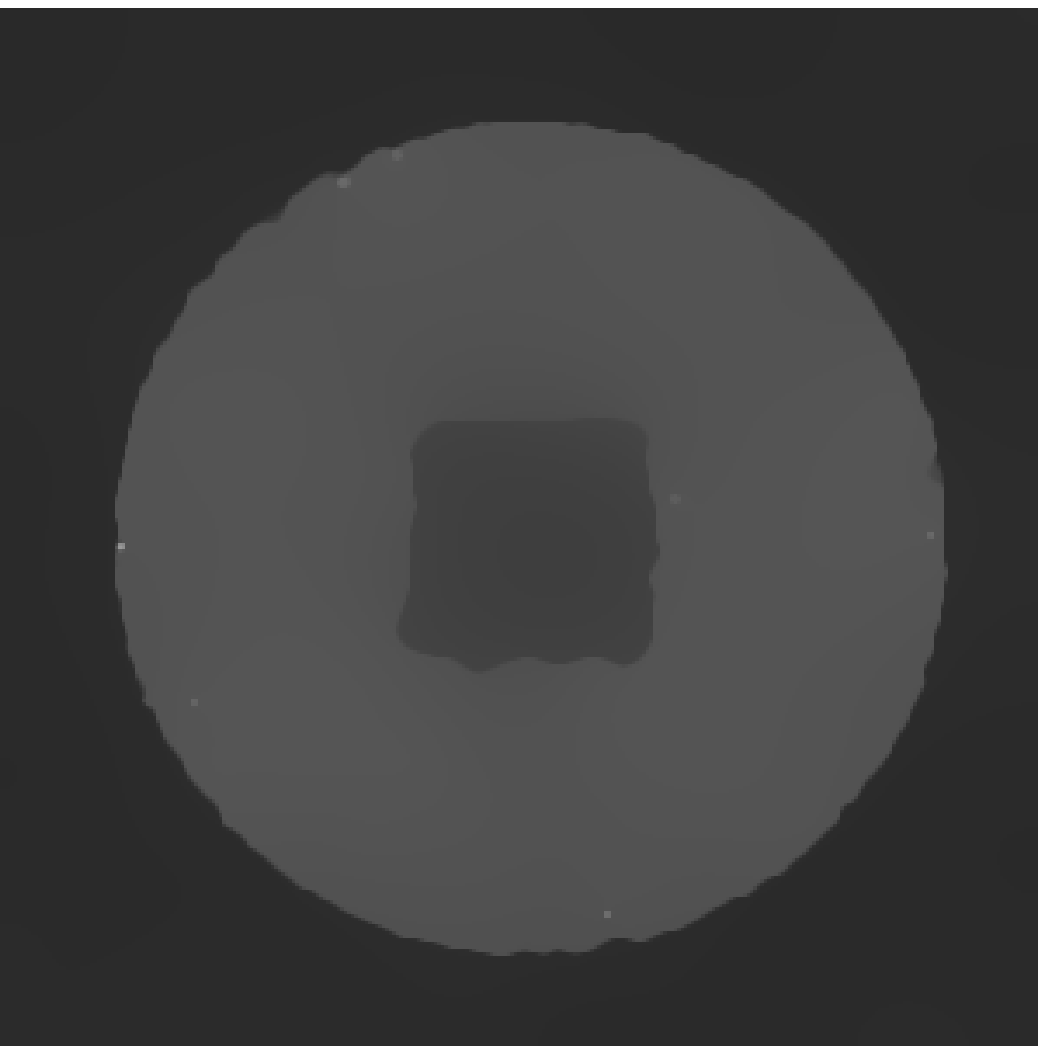}               
                \caption{Proposed}
                \label{fig:circle_3_tdm}
       \end{subfigure}
     
\caption{Image corrupted with speckle look L=3 and restored by different models. }\label{circle_3}
\end{figure}
%======================================================================================================
%======================================================================================================
\begin{figure}
       \centering
          \begin{subfigure}[b]{0.25\textwidth}           
                \includegraphics[scale=0.37]{circle}               
                \caption{Original}
                \label{fig3:circle}
       \end{subfigure}%
       \begin{subfigure}[b]{0.25\textwidth}           
                \includegraphics[scale=0.37]{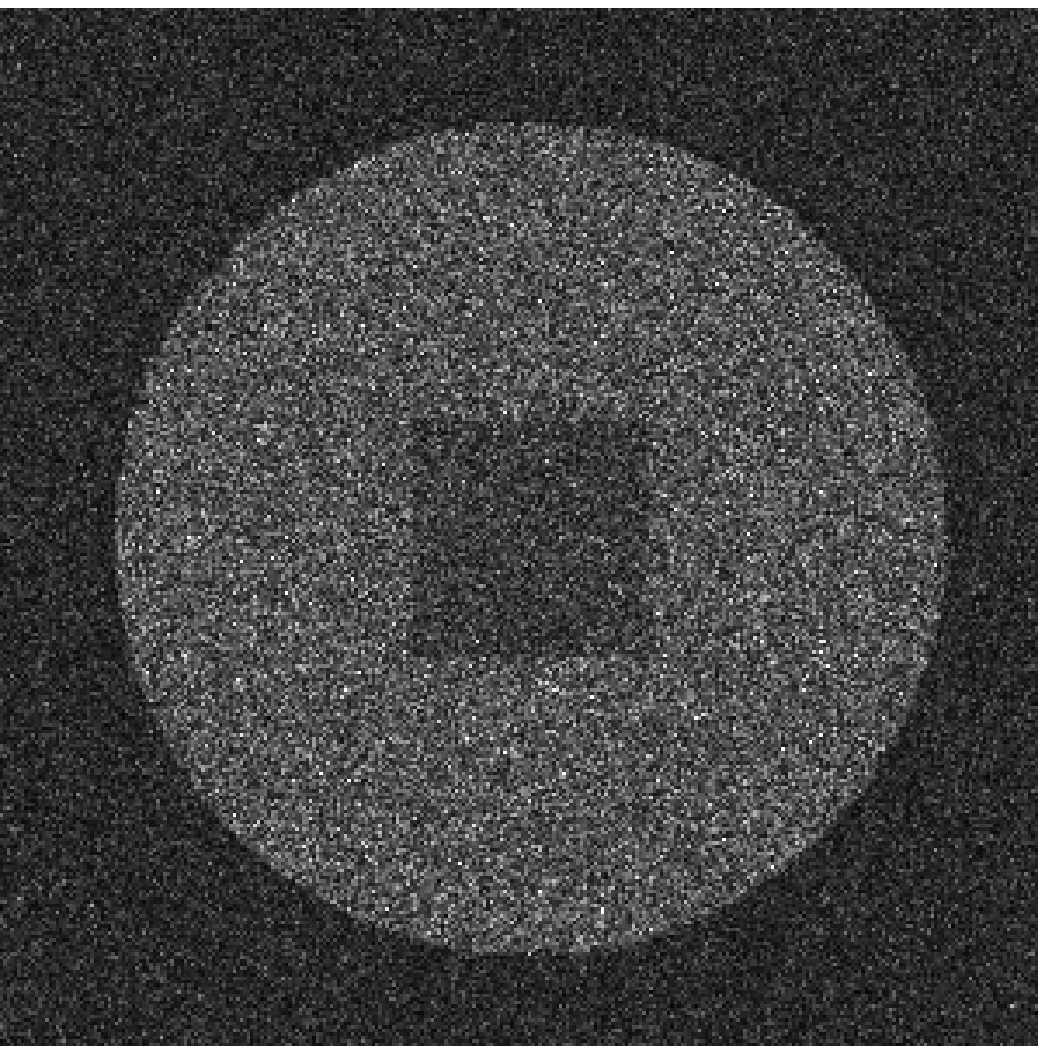}               
                \caption{Noisy}
                \label{fig:circle_5}
       \end{subfigure}% 
         \begin{subfigure}[b]{0.25\textwidth}           
                \includegraphics[scale=0.37]{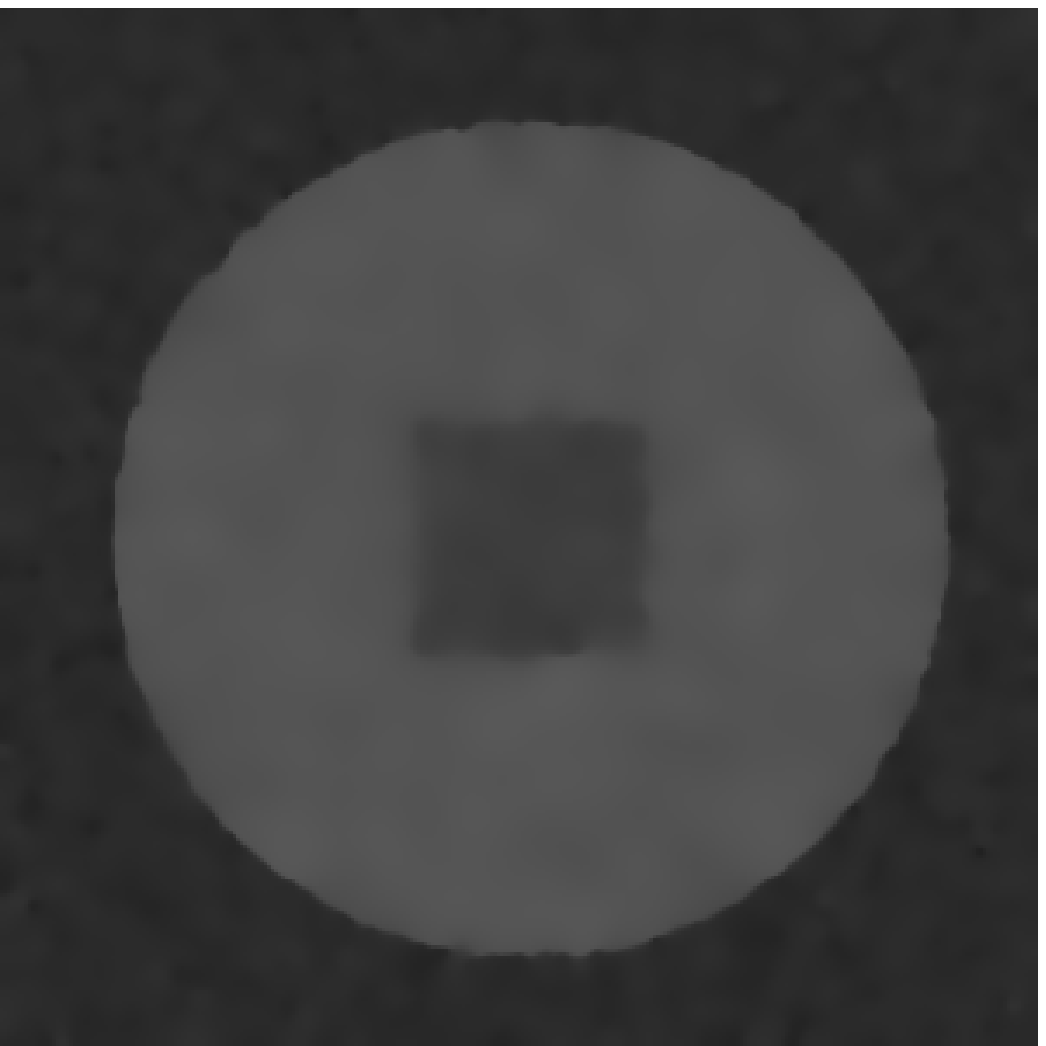}               
                \caption{Shan}
                \label{fig:circle_5_zzdb}
       \end{subfigure}% 
        \begin{subfigure}[b]{0.25\textwidth}           
                \includegraphics[scale=0.37]{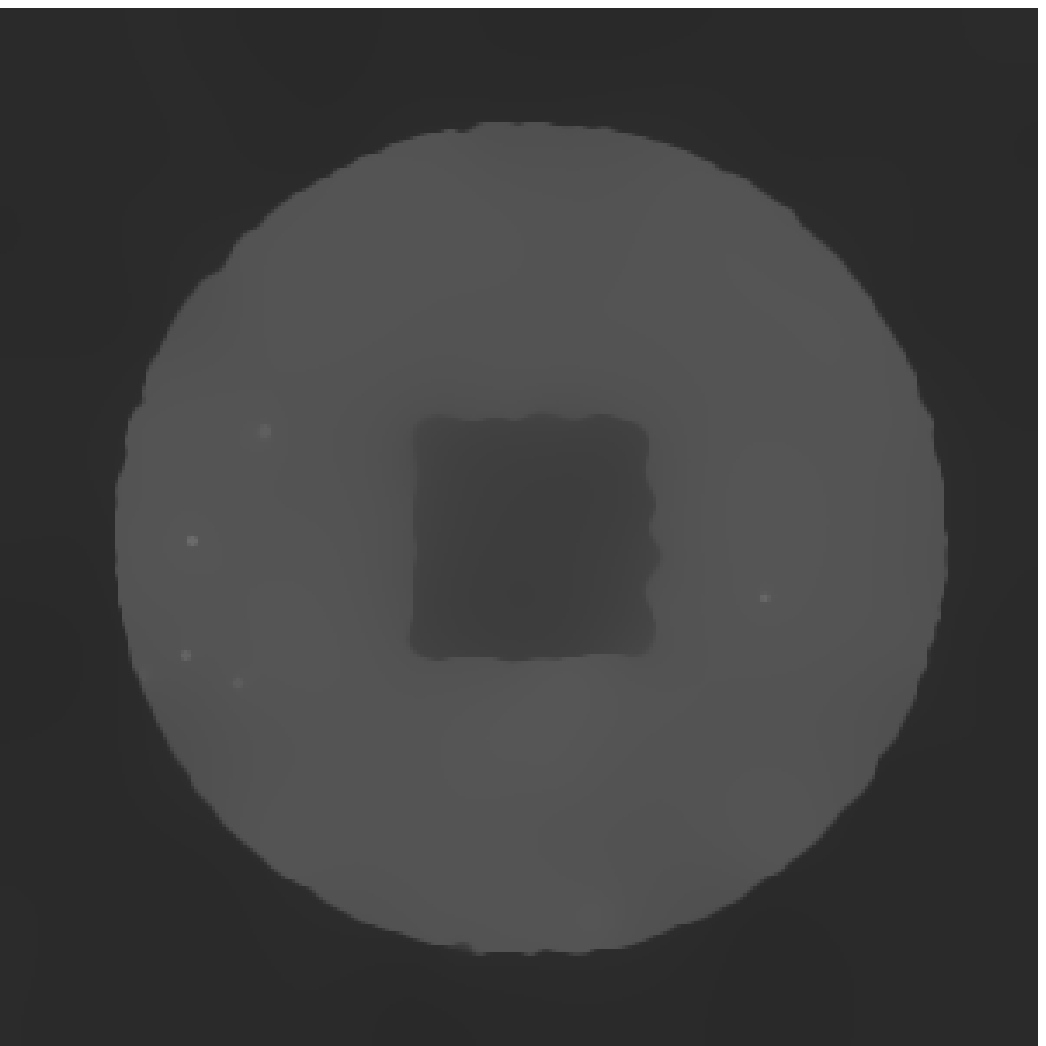}               
                \caption{Proposed}
                \label{fig:circle_5_tdm}
       \end{subfigure}
     
\caption{Image corrupted with speckle look L=5 and restored by different models. }\label{circle_5}
\end{figure}
%==============================================================================================================
\begin{figure}
       \centering
        \begin{subfigure}[b]{0.25\textwidth}           
                \includegraphics[scale=0.37]{circle}               
                \caption{Original}
                \label{fig4:circle}
       \end{subfigure}%
       \begin{subfigure}[b]{0.25\textwidth}           
                \includegraphics[scale=0.37]{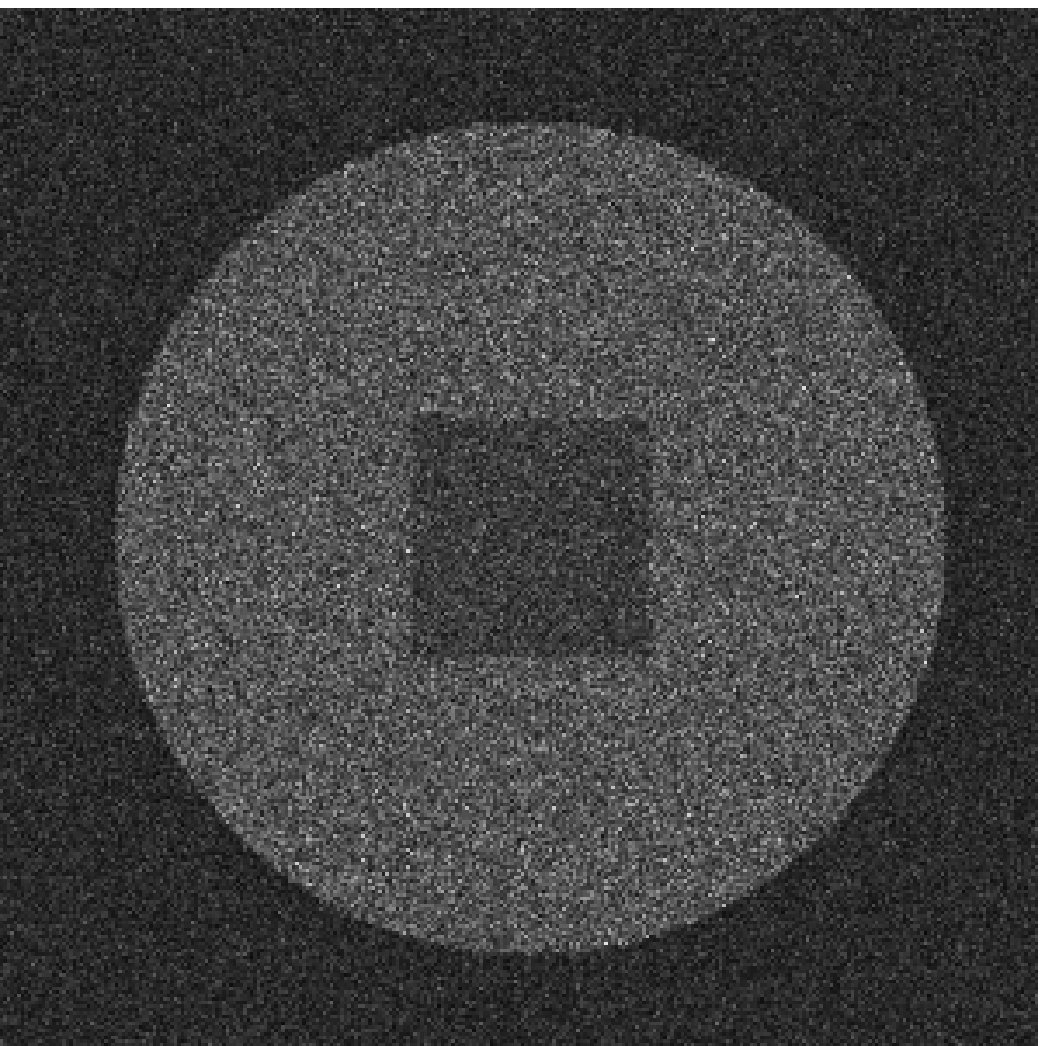}               
                \caption{Noisy}
                \label{fig:circle_10}
       \end{subfigure}% 
         \begin{subfigure}[b]{0.25\textwidth}           
                \includegraphics[scale=0.37]{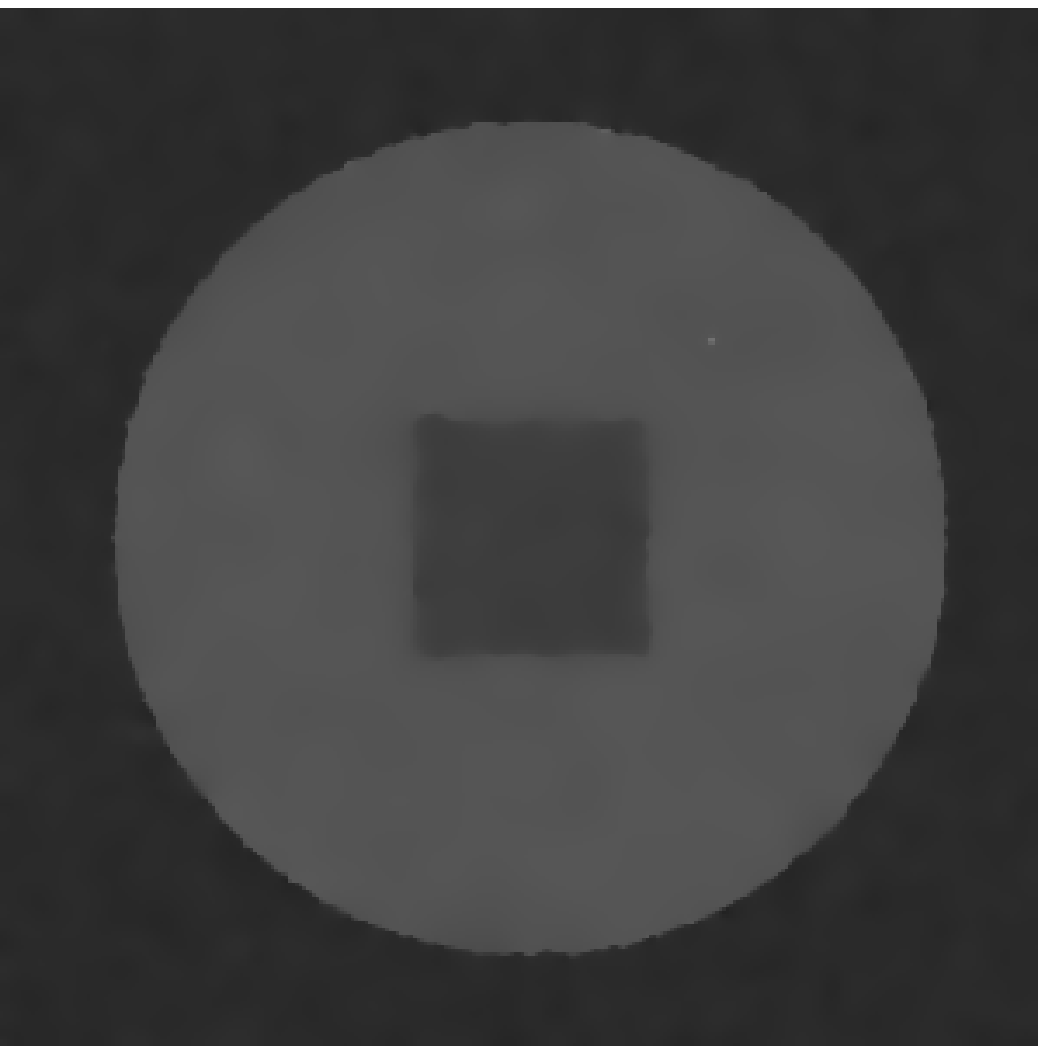}               
                \caption{Shan}
                \label{fig:circle_10_zzdb}
       \end{subfigure}% 
        \begin{subfigure}[b]{0.25\textwidth}           
                \includegraphics[scale=0.37]{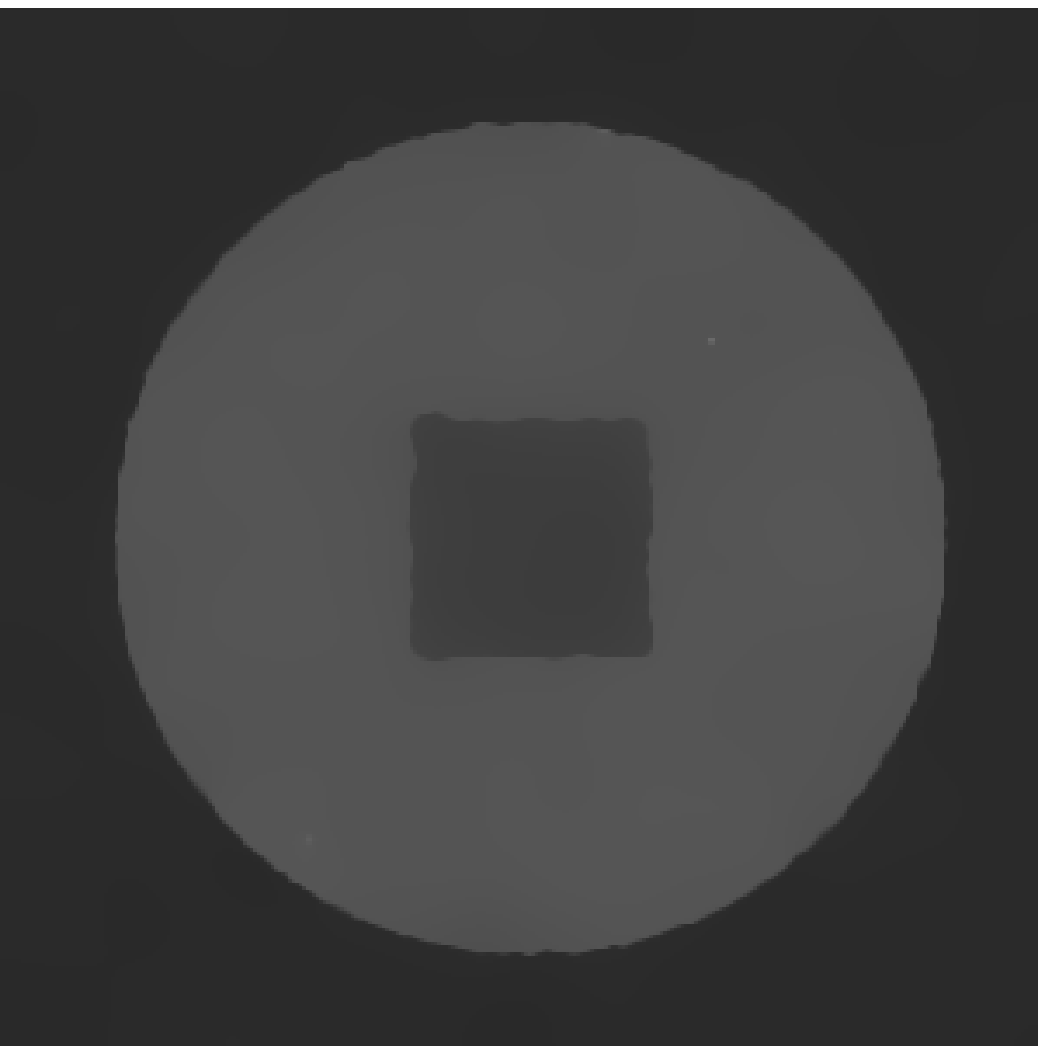}               
                \caption{Proposed}
                \label{fig:circle_10_tdm}
       \end{subfigure}
     
%\caption{Image corrupted with speckle look L=10 and restored by different models. }\label{fig:circle_10}
%\end{figure}
%==============================================================================================================
%\begin{figure}
%       \centering
       \begin{subfigure}[b]{0.35\textwidth}           
                \includegraphics[scale=0.45]{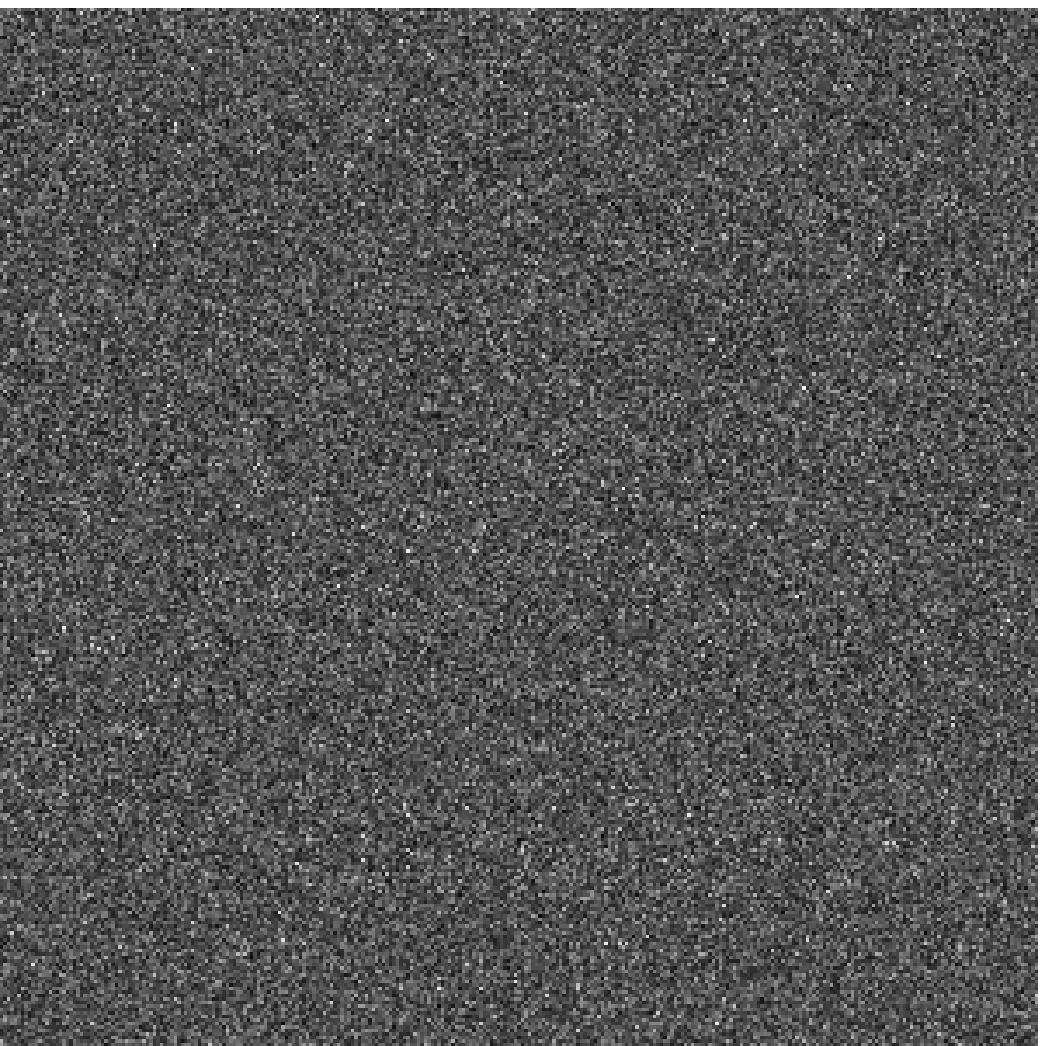}               
                \caption{Original}
                \label{fig:circle_ratio}
       \end{subfigure}% 
         \begin{subfigure}[b]{0.35\textwidth}           
                \includegraphics[scale=0.45]{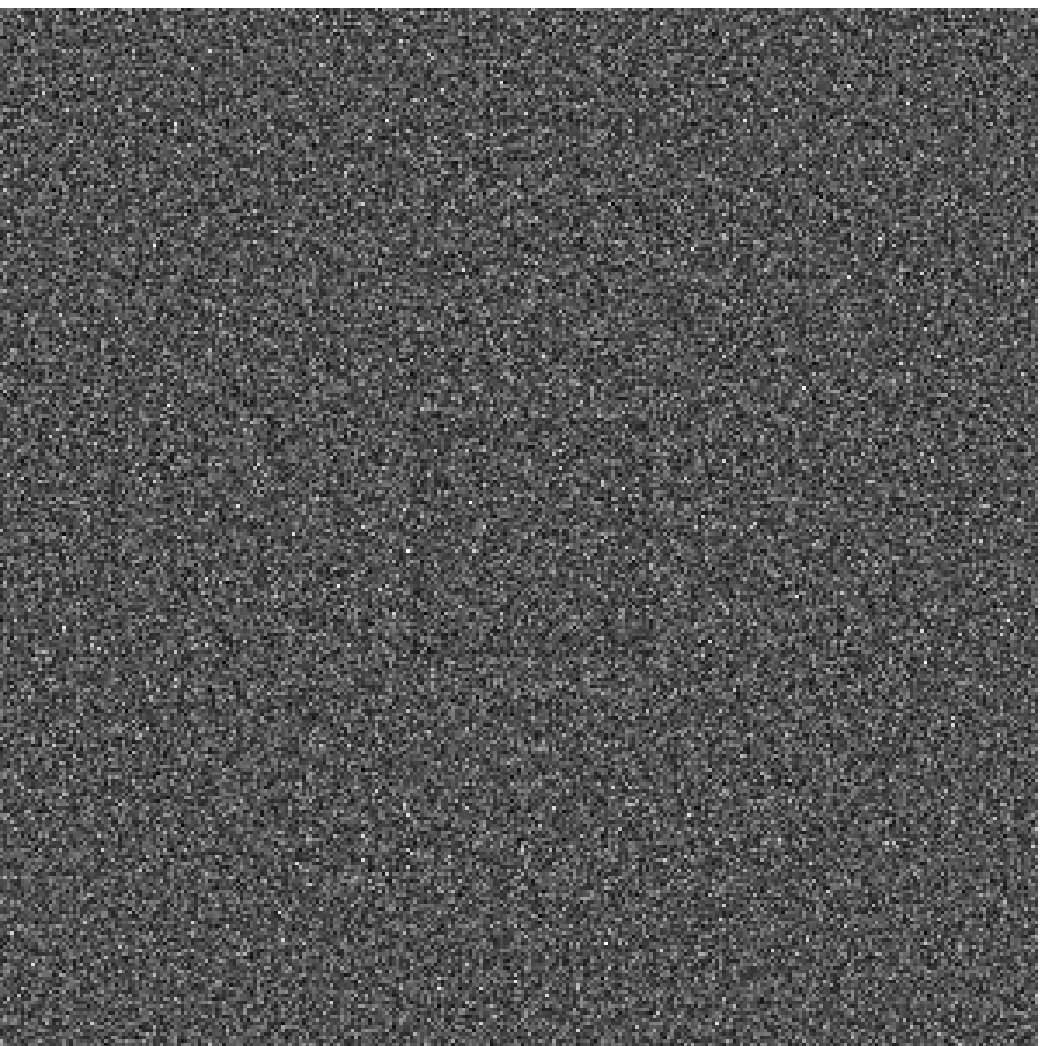}               
                \caption{Shan}
                \label{fig:circle_10_zzdb_ratio}
       \end{subfigure}% 
        \begin{subfigure}[b]{0.35\textwidth}           
                \includegraphics[scale=0.45]{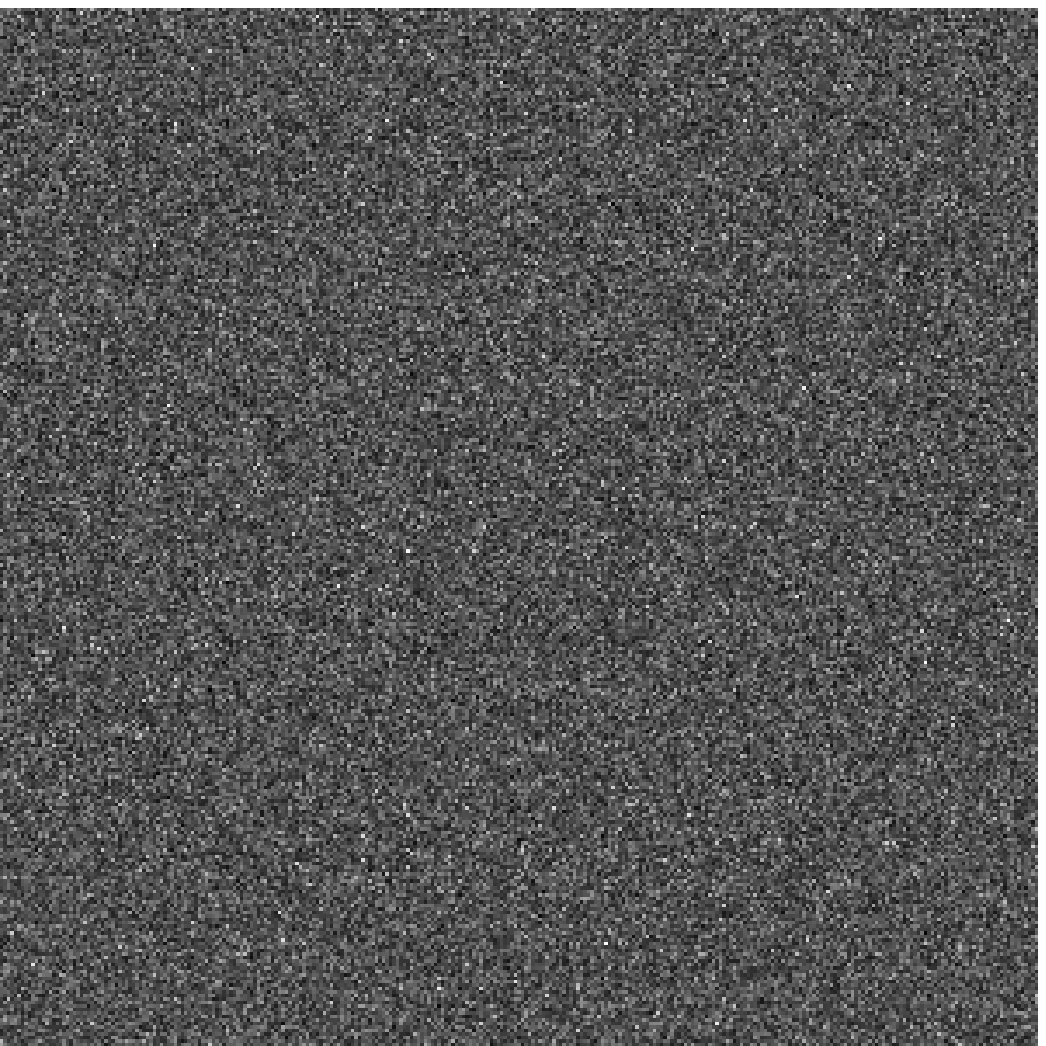}               
                \caption{Proposed}
                \label{fig:circle_10_tdm_ratio}
       \end{subfigure}
     
\caption{Upper row: Image corrupted with speckle look L=10 and restored images. Lower row: Ratio images. }\label{fig:circle_10_ratio}
\end{figure}
%===========================================================================================================
\begin{figure}
       \centering
       
           \begin{subfigure}[b]{0.4\textwidth}           
           \includegraphics[scale=0.45]{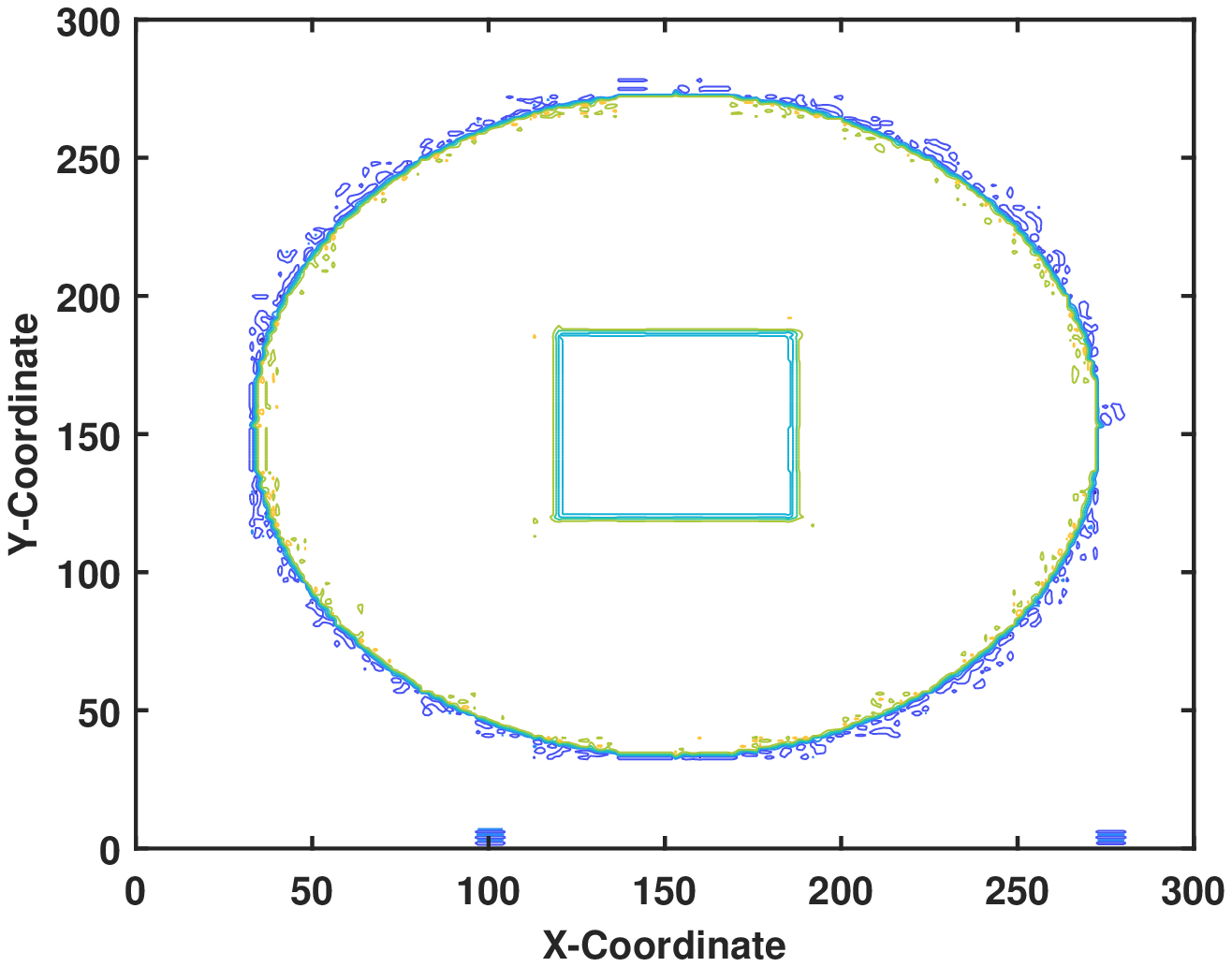}           
                \caption{Original}
                \label{fig:3a_cont}
        \end{subfigure}%
              \begin{subfigure}[b]{0.4\textwidth}           
                \includegraphics[scale=0.45]{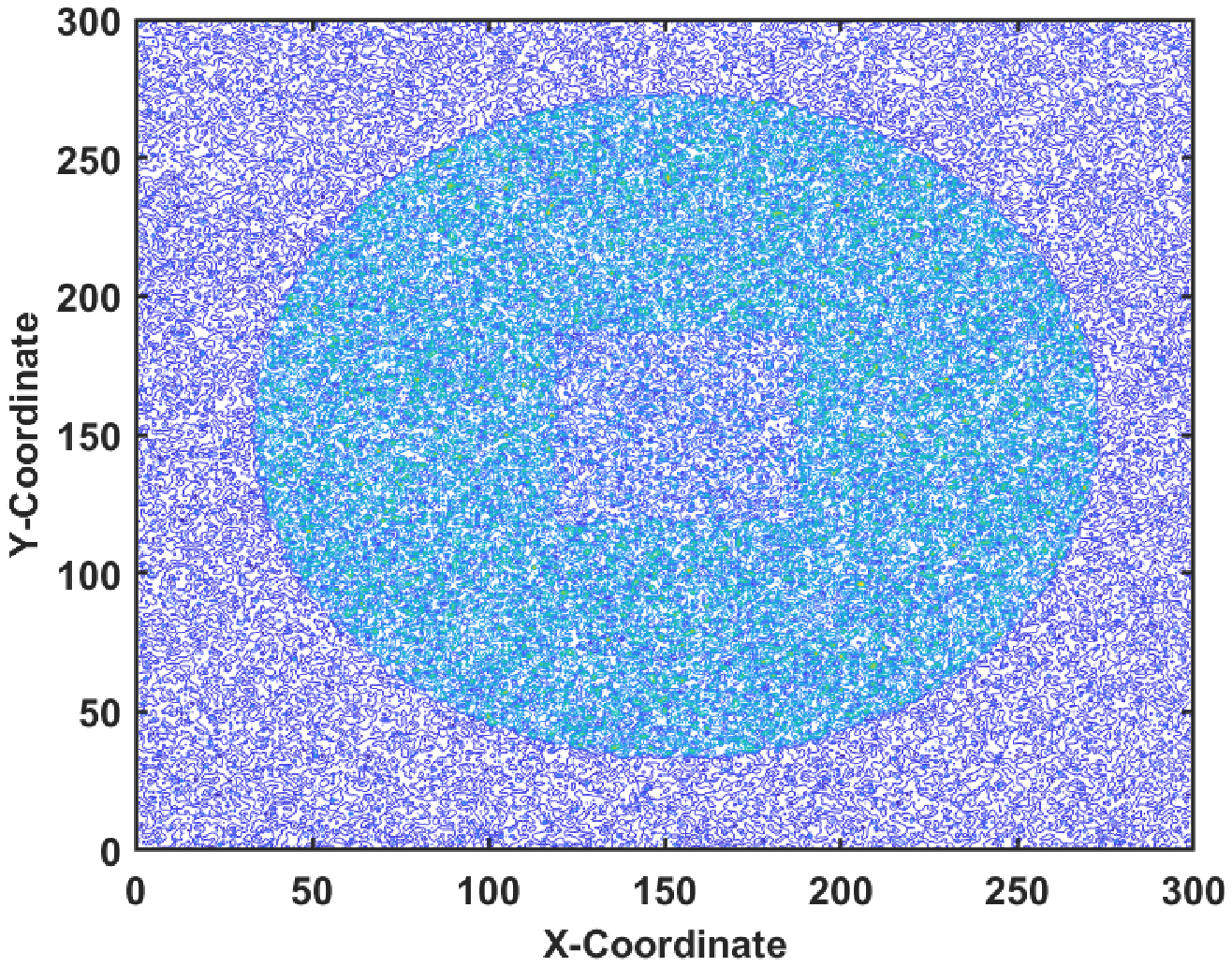}               
                \caption{Noisy}
                \label{fig:3b_cont}
       \end{subfigure}% 
       
        \begin{subfigure}[b]{0.4\textwidth}           
                \includegraphics[scale=0.45]{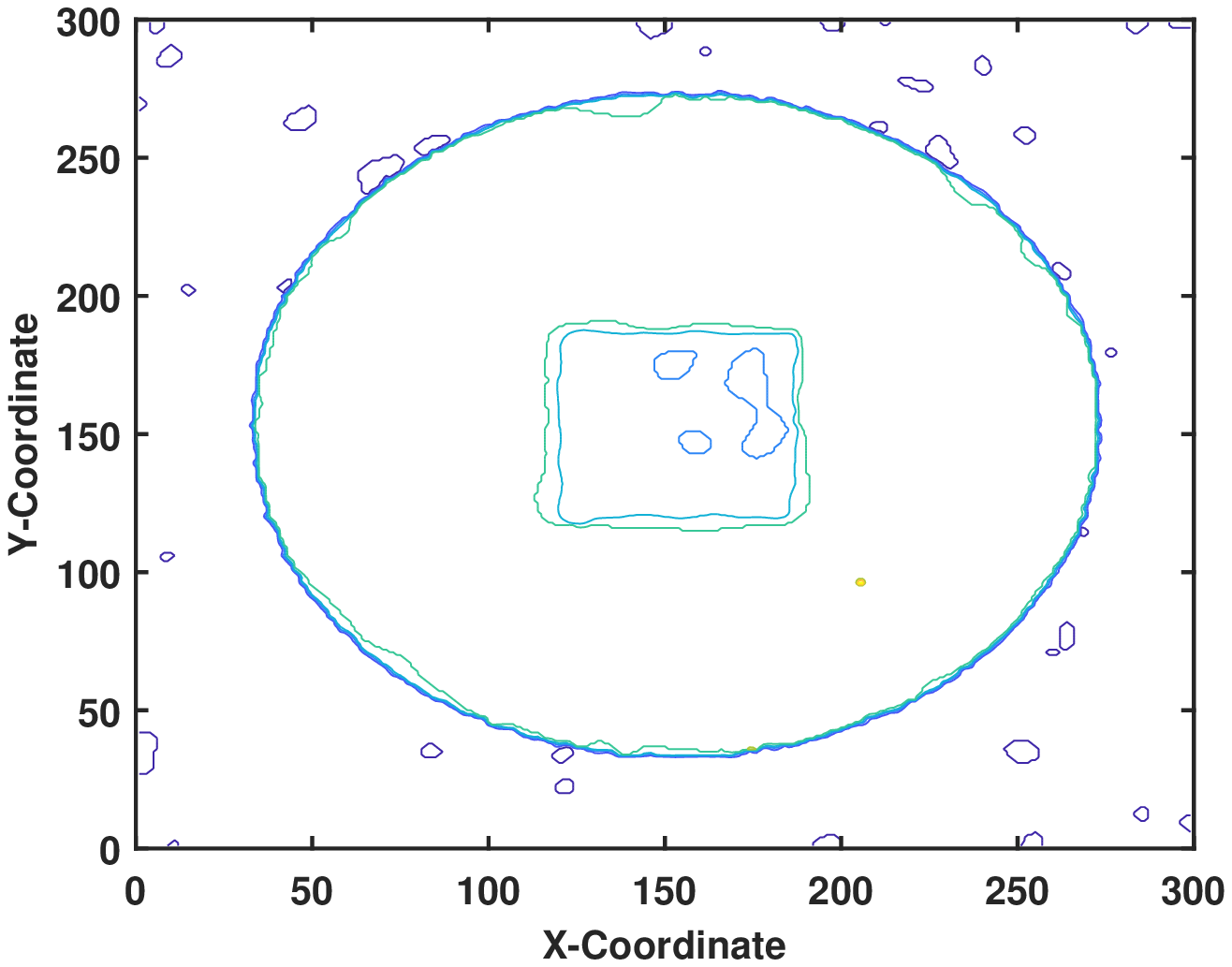}               
                \caption{Shan}
                \label{fig:3d_cont}
         \end{subfigure}% 
        \begin{subfigure}[b]{0.4\textwidth}           
                \includegraphics[scale=0.45]{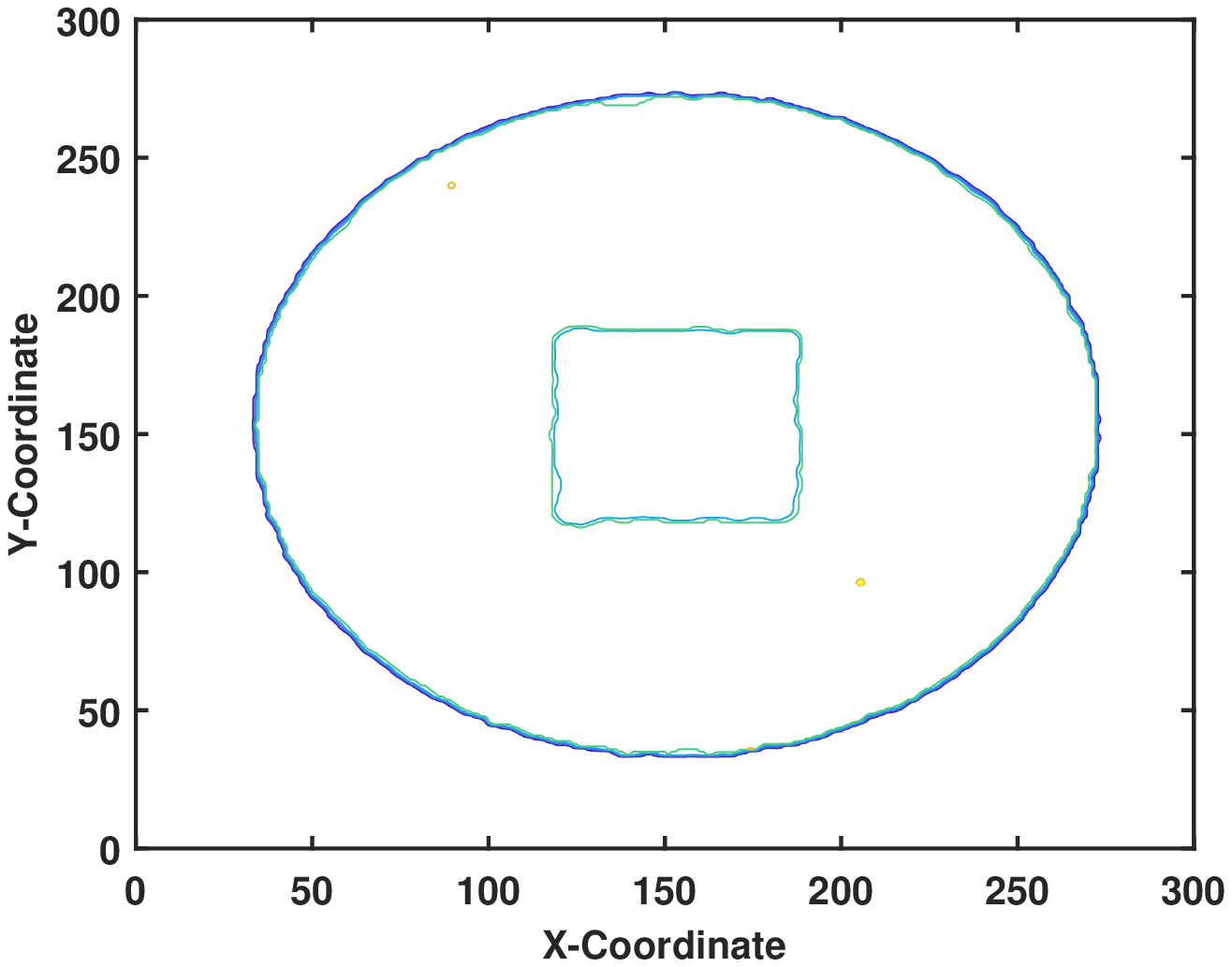}               
                \caption{Proposed}
                \label{fig:3e_cont}
        \end{subfigure}% 
       
\caption{Contour maps of the restored imgaes in figure \ref{fig:circle_10_ratio}.}\label{fig:circle_10_cont}
\end{figure}
%======================================================================================================
\begin{figure}
       \centering
       
           \begin{subfigure}[b]{0.4\textwidth}           
           \includegraphics[scale=0.5]{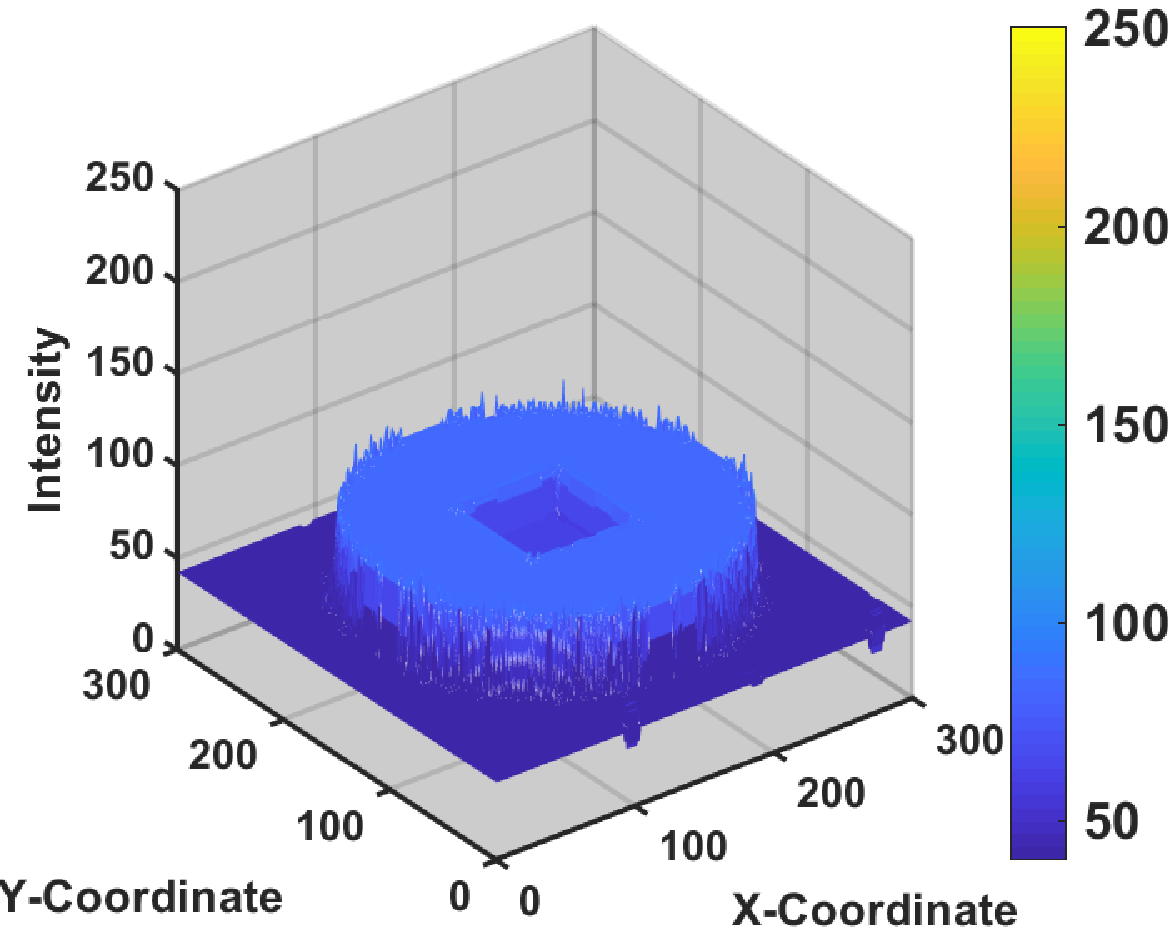}           
                \caption{Original}
                \label{fig:2a_3d}
        \end{subfigure}%
              \begin{subfigure}[b]{0.4\textwidth}           
                \includegraphics[scale=0.5]{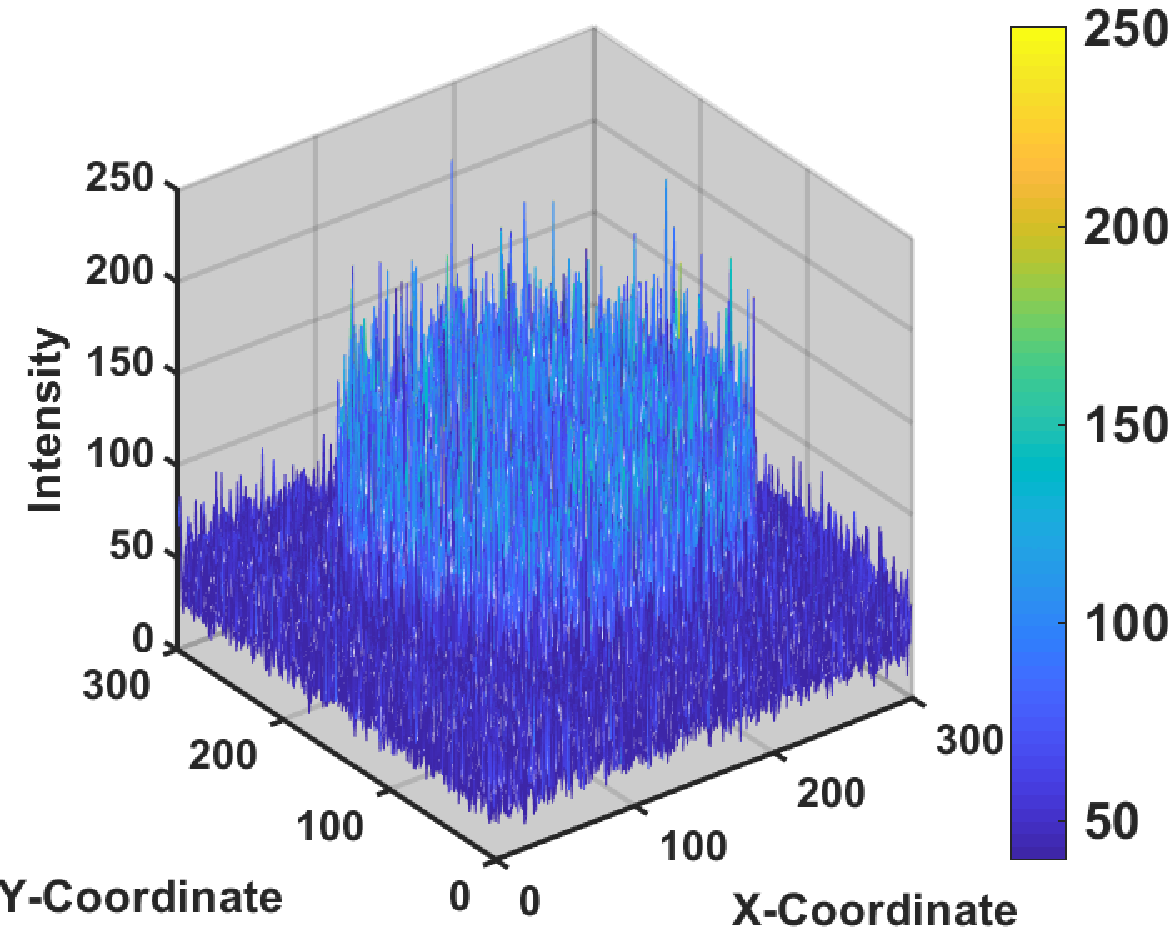}               
                \caption{Noisy}
                \label{fig:2b_3d}
       \end{subfigure}% 

        \begin{subfigure}[b]{0.4\textwidth}           
                \includegraphics[scale=0.5]{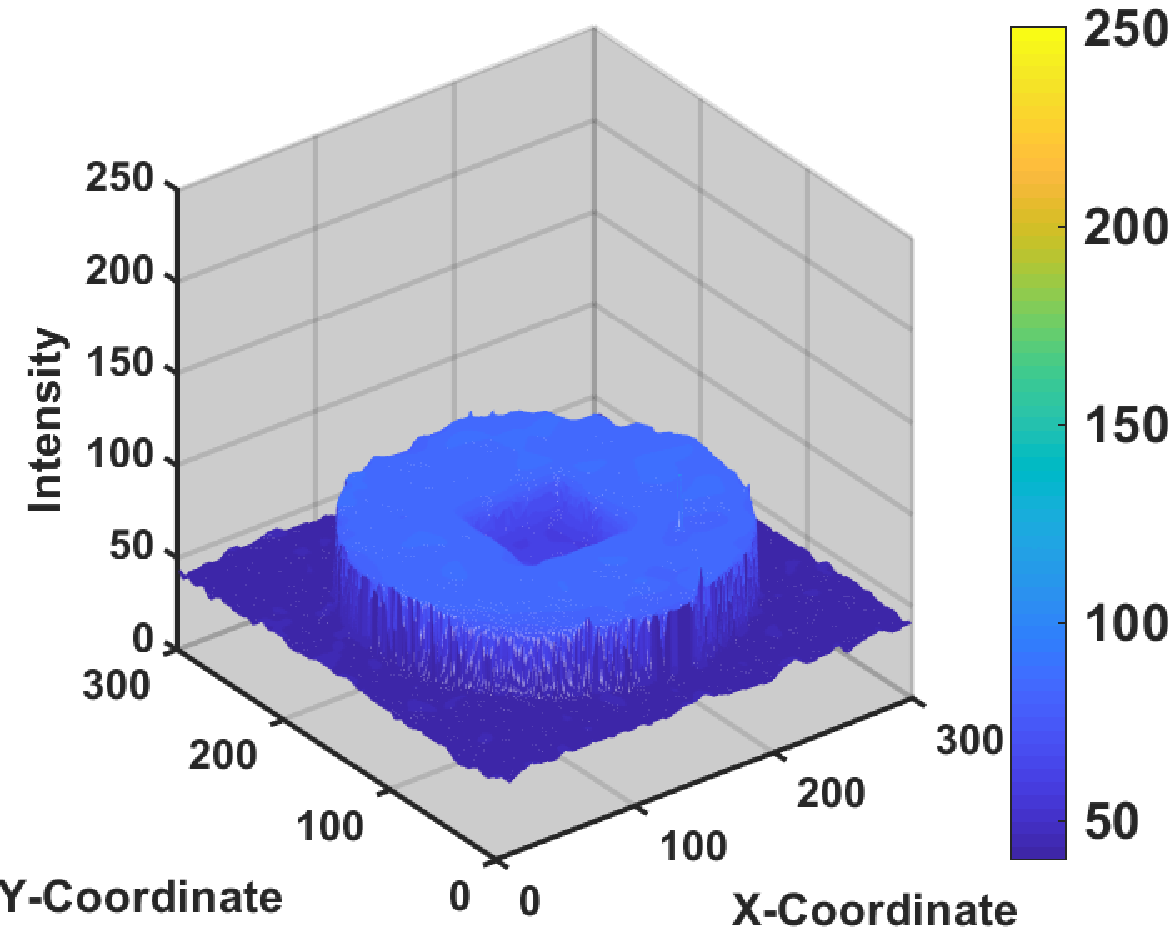}               
                \caption{Shan}
                \label{fig:2d_3d}
         \end{subfigure}% 
        \begin{subfigure}[b]{0.4\textwidth}           
                \includegraphics[scale=0.5]{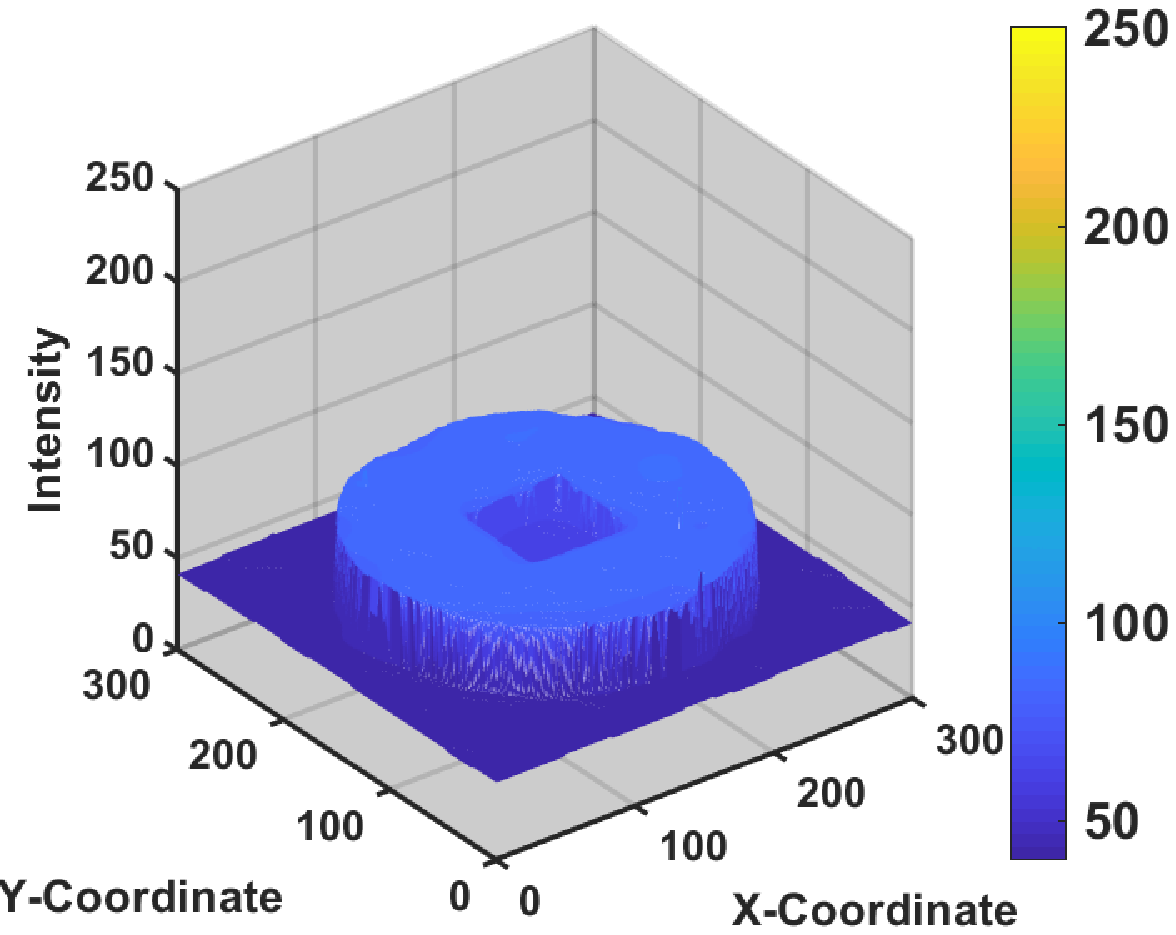}               
                \caption{Proposed}
                \label{fig:2e_3d}
        \end{subfigure}% 
       
\caption{ 3D surface plots of the restored imgaes in figure \ref{fig:circle_10_ratio}.}\label{fig:circle_10_3d}
\end{figure}
%\fi
%================================================================================================================
\begin{center}
\begin{table}
\caption{Left table: Comparison of SSIM and PSNR values of despeckled images. Right table: Parameter values for the numerical experiments.}
\label{tab:psnr_ssim_parameter}
%\centering
\begin{tabular}{ll}
\scalebox{0.8}{
%\begin{tabular}{ccrrrrrrrrrrrrrr}
\begin{tabular}[t]{llccccccccc}
\toprule
    \multirow{2}[8]{*}{Image} & \multirow{2}[8]{*}{$L$}  & \multicolumn{2}{c}{Shan Model\cite{shan2019multiplicative}}  & \multicolumn{2}{c}{Proposed Model} \\
\cmidrule(r){3-4}
		\cmidrule(r){5-6}
	%	\cmidrule(r){7-8}
	%	\cmidrule(r){9-10}  
	%	\cmidrule(r){11-12} 

&         & \multicolumn{1}{c}{SSIM} & \multicolumn{1}{c}{PSNR} & \multicolumn{1}{c}{SSIM}  & \multicolumn{1}{c}{PSNR}  \\
    \midrule
Boat    & 1             & 0.5975 & 17.10  & \textbf{0.6096} & \textbf{17.12} \\
        & 3             & 0.7087 & 22.71  & \textbf{0.7347} & \textbf{22.85} \\
        & 5             & 0.7508 & 24.73  & \textbf{0.7905} & \textbf{24.93} \\
        & 10            & 0.8325 & 26.98  & \textbf{0.8422} & \textbf{27.12} \\
        & 33            & 0.8941 & 29.57  & \textbf{0.9057} & \textbf{29.72} \\
        &               &        &       &        &        &                 &        \\
Brick   & 1             & 0.2930 & 12.17  & \textbf{0.2954} & \textbf{12.19} \\
        & 3             & 0.3837 & 17.08  & \textbf{0.3861} & \textbf{17.11} \\
        & 5             & 0.4291 & 19.34  & \textbf{0.4355} & \textbf{19.40} \\
        & 10            & 0.4947 & 22.06  & \textbf{0.4960} & \textbf{22.18} \\
        & 33            & 0.5943 & 25.40  & \textbf{0.5961} & \textbf{25.53} \\
        &               &        &       &        &        &                 &       \\
Circle  & 1             & 0.9582 & 34.30  & \textbf{0.9644} & \textbf{34.70} \\
        & 3             & 0.9735 & 38.10  & \textbf{0.9772} & \textbf{39.53} \\
        & 5             & 0.9765 & 39.36  & \textbf{0.9806} & \textbf{40.73} \\
        & 10            & 0.9817 & 41.26  & \textbf{0.9865} & \textbf{42.85} \\
        & 33            & 0.9870 & 43.64  & \textbf{0.9889} & \textbf{44.62} \\
        &               &        &        &                 &                \\
\midrule        
%\bottomrule
\end{tabular}
}
&
\scalebox{0.87}{
\begin{tabular}[t]{llcccccc}
\toprule
    \multirow{2}[4]{*}{Image} & \multirow{2}[4]{*}{ $L$ }  & \multicolumn{2}{c}{Shan\cite{shan2019multiplicative}} & \multicolumn{3}{c}{Proposed} \\
        \cmidrule(r){3-4}
		\cmidrule(r){5-7}

		 &  &  \multicolumn{1}{c}{$\alpha$}  & \multicolumn{1}{c}{$\beta$}  &\multicolumn{1}{c}{$\gamma$} & \multicolumn{1}{c}{$\nu$} & \multicolumn{1}{c}{$K$}   \\
  \midrule
Boat        & 1                  & 1    &  1   & 5  & 1   &  2        \\
            & 3                  & 1.2  &  1   & 4  & 1.5 &  2        \\
            & 5                  & 1.3  &  1   & 2  & 1.5 &  1        \\
            & 10                 & 1.4  &  1.2 & 2  & 2   &  1       \\
            & 33                 & 1.5  &  1.5 & 2  & 3   &  1      \\
\midrule
Brick       & 1                 & 1    & 1    & 5  & 1   &  4        \\
            & 3                 & 1.2  & 1    & 4  & 1.3 &  3        \\
            & 5                 & 1.4  & 1    & 2  & 1.5 &  2        \\
            & 10                & 1.6  & 1    & 2  & 2   &  1        \\
            & 33                & 1.7  & 1    & 2  & 3   &  1      \\
\midrule
Circle      & 1                 & 1.5  & 2    & 10 & 1   &  1       \\
            & 3                 & 1.5  & 2    & 10 & 1   &  1        \\
            & 5                 & 2    & 2.25 & 5  & 1   &  1        \\
            & 10                & 2    & 2.25 & 2  & 1   &  1        \\
            & 33                & 2    & 2.5  & 2  & 1   &  1      \\

\bottomrule
\end{tabular}
}
\end{tabular}
\end{table}
\end{center}
%====================================================================================================================
\section{Conclusion}
\label{sec:Conclusion}

This work suggests an efficient telegraph diffusion-based multiplicative speckle noise removal model. Such a new method intends to preserve the image edges during the noise removal process. To overcome the limitations of gradient-based despeckling models as well as parabolic PDE based models, we considered a hybrid approach. Here we combine a gray level indicator function with gradient-based diffusion in a telegraph diffusion framework for image restoration. To the best of our knowledge, the gray level indicator based telegraph diffusion model has not been used before for speckle noise suppression.
Also, we established the existence and uniqueness of a weak solution to the suggested model using Schauder's fixed point theorem. Moreover, we prove the boundedness of the weak solution.
Numerical experiments have been conducted to highlight the efficiency of the proposed model for despeckling using different types of test images.  Computational result of the present model compares with a recently developed model. From the experiment results of the proposed model, we can conclude that the images are suitably recovered without introducing undesired artifacts.  A potential direction that the telegraph diffusion model can be extended to handle texture preservation issues in various real-life images, which are degraded by mixed noises. Another significant step might be the study of the advanced numerical solver to enhance the convergence speed of the proposed model.

\thispagestyle{empty}

\bibliographystyle{unsrt}

\end{document}